\documentclass[10pt,a4paper,draft]{article}
\usepackage{amsmath}
\usepackage{amsfonts}
\usepackage{amssymb}
\usepackage{a4,a4wide}
\usepackage{graphics}
\usepackage{color}
\usepackage[T1]{fontenc}

\def\ep{{\varepsilon}}
\def\R{\mathbb R}
\def\mass{m_\Gamma}

\newcommand{\fdem}{ \hfill $\square$ }

\newtheorem{theo}{\textbf{Theorem}}[section]

\newtheorem{lem}[theo]{\textbf{Lemma}}
\newtheorem{prop}[theo]{\textbf{Proposition}}

\newtheorem{defi}[theo]{\textbf{Definition}}

\newtheorem{ass}[theo]{\textbf{Assumption}}

\newtheorem{rem}[theo]{\textbf{Remark}}

\def\e{{\mathcal{E}}}

\title{Travelling waves in a nonlocal reaction-diffusion equation as a model for a population structured by
a space variable and a phenotypical trait}
\date{}

\begin{document}

\maketitle

\begin{center}
{\large\bf Matthieu Alfaro \footnote{ I3M, Universit\'e de
Montpellier 2, CC051, Place Eug\`ene Bataillon, 34095 Montpellier
Cedex 5, France. E-mail: malfaro@math.univ-montp2.fr},
J{\'e}r{\^o}me Coville \footnote{Equipe BIOSP, INRA Avignon,
Domaine Saint Paul, Site Agroparc, 84914 Avignon Cedex 9, France.
E-mail: jerome.coville@avignon.inra.fr} and  Ga\"el Raoul
\footnote{Centre d'\'Ecologie Fonctionnelle et \'Evolutive, UMR
5175, CNRS, 1919 Route de Mende, 34293 Montpellier, France.
E-mail: raoul@cefe.cnrs.fr
.}.}\\
[2ex]

\end{center}



\vspace{10pt}

\begin{abstract} We consider a nonlocal reaction-diffusion equation as a model for a population structured by a space
variable and a phenotypical trait. To sustain the possibility of
invasion in the case where an underlying principal eigenvalue is
negative, we investigate the existence of travelling wave
solutions. We identify a minimal speed $c^*>0$, and prove the
existence of waves when $c\geq c^*$ and the non existence when
$0\leq c<c^*$.\\


\end{abstract}

\section{Introduction}\label{s:intro}

\subsection{Setting of the problem}

In this paper we are interested in propagation phenomena for
nonlocal reaction-diffusion equations of the form
\begin{equation}\label{eq0}
\partial_t n(t,x,y)-\Delta _{x,y}n(t,x,y)=\left(r(y-Bx\cdot e)-\int _{\R}k(y-Bx\cdot e,y'-Bx\cdot e) n(t,x,y')\,dy'\right)n(t,x,y),
\end{equation}
where $(x,y)\in \R^{d}\times \R$, $e\in S^{d-1}$, $B\geq 0$, $r:
\R\to \R$ and $k:\R^2 \to \R ^+$.

Such equations have appeared in some population dynamic models,
see \cite{Peck}, \cite{Prevost}, \cite{Polechova}, \cite{Mir-Rao}.
In this context, $n(t,x,y)$ denotes a density of population
structured by a spatial variable $x\in \R^d$ and by a phenotypical
trait $y\in\R$.
 The population  is then  submitted to four essential processes:  spatial dispersion, mutations, growth and competition. The spatial dispersion and the mutations are modelled by
diffusion operators.  The growth rate  of the population at a
location $x$ and trait $y$ is given --- for all times $t$--- by
$r(y-Bx\cdot e)$, where $r$ is typically negative outside a
bounded interval. This corresponds to a population living in an
environmental cline: to survive at the location $x$, an individual
must have a trait close to $y=Bx\cdot e$. Therefore, to be able to
invade the environment, the population needs to evolve. Finally,
we consider a logistic regulation
 of the population density that is local in the spatial variable  and nonlocal in  the trait. In other words,  we consider that there
 exists an intra-specific competition (for e.g. food) at each location, which may depend on the traits of the competitors. For a
 rigorous derivation of this model from individual based models, we refer to
 \cite{Champagnat}. In Section \ref{intro_bio}, we will discuss in
 more details the biological aspects of our work.

The existence of  global  solutions for the Cauchy problem and  of
non trivial steady states for \eqref{eq0} have been investigated
respectively in \cite{Prevost} and  \cite{Arnold}. Also, numerical
simulations (see e.g. \cite{Peck}, \cite{Prevost},
\cite{Polechova}) show that the population can either
 go extinct, or propagate while adapting to local environments.

The aim of this work is to analyze this propagation phenomena
through the study of travelling front solutions. The travelling
front solutions are particular  solution of \eqref{eq0} describing
the transition at a  constant speed $c$ from one stationary
solution to another one. Such solutions have  proved in numerous
situations their utility in describing the dynamics of a
population modelled by a reaction diffusion equation.  In the case
of the classical Fisher-KPP equation
\begin{equation}
 \partial_tn-\Delta n=\left(1- n\right)n, \label{eq-kpp}
\end{equation}
we refer among others to \cite{Fis}, \cite{KPP}, \cite{Aro-Wei2}
\cite{W1}: there exists planar fronts $\phi(x.e-ct)$ connecting
$0$ to $1$, for all speed $c\ge c^*= \sqrt{2}$. Moreover, the
minimal speed of the front $c^*$ corresponds to the so called {\it
spreading speed} of propagation. Travelling front solutions in
heterogeneous versions of \eqref{eq-kpp} with periodicity in
space, in time, or more general media are studied  in
\cite{berestycki-nirenberg}, \cite{HZ}, \cite{Xin3},
\cite{berestycki-hamel-cpam}, \cite{Berestycki-Hamel-Roques-II},
\cite{nadin}, \cite{Nolen-Ryzhik}. Nonlocal versions of
\eqref{eq-kpp} where the Laplace operator is replaced by a
nonlocal operator are studied in \cite{CD1}, \cite{CDM2, CDM3},
\cite{Shen2012}. For very general reaction diffusion equations, we
refer to \cite{berestycki-hamel-cpam2} for a definition of
generalized transition waves and their properties.

  It is
worth noticing that when the competition term is replaced by a
local (in $x$ and $y$) density regulation, equation \eqref{eq0}
becomes the following heterogeneous reaction diffusion equation
\begin{equation}
\partial_tn(t,x,y)-\Delta _{x,y}n(t,x,y)=\left(r(y-Bx\cdot e)-h(y-Bx\cdot e) n(t,x,y)\right)n(t,x,y), \label{eq-ber-cha}
\end{equation}
which was recently investigated by Berestycki and Chapuisat
\cite{Ber-Cha}: they prove the existence of a critical speed
$c^*$, for which there exists a  travelling front of
\eqref{eq-ber-cha} for any  speed $c\ge c^*$.

As far as nonlocal equations of the form \eqref{eq0} are
concerned, far less seems to be known in the literature. The
travelling wave analysis has been done either using a formal
Hamilton-Jacobi approach \cite{Bouin} for a model close to
\eqref{eq0}, or for a population structured by one variable only,
that is $n=n(t,x)$, submitted to a nonlocal competition
\cite{Ber-Nad-Per-Ryz}. To our knowledge, there is no  result on
the existence of travelling waves for \eqref{eq0} and related
models.

\subsection{Assumptions and main results}

As suggested by the numerical simulations mentioned above, we
expect that during an invasion, the
 population adapts locally to the environmental gradient. To observe travelling waves of \eqref{eq0}, we therefore perform the change of variable
\begin{equation}
\tilde n(t,x,z)=n(t,x,z+Bx\cdot e).
\end{equation}
Then \eqref{eq0} is recast as
\begin{equation}\label{eq1}
\partial_t \tilde n(t,x,z)-\tilde{\mathcal E}(\tilde n)(t,x,z)=\left(r(z)-\int _{\R} k(z,z')\tilde n(t,x,z')\,dz'\right)\tilde
n(t,x,z),
\end{equation}
where $\tilde {\mathcal E}(\tilde n):=\Delta_x\tilde
n+\left(B^2+1\right)\tilde n_{zz}-2B\partial_z(\nabla_x \tilde
n\cdot e)$.  Since \eqref{eq0} is invariant under any rotation in
$\R^d$, without loss of generality we can assume that $e=e_1$
 so that $\tilde{\e}(\tilde n)= \Delta_x\tilde n+\left(B^2+1\right)\tilde
n_{zz}-2B\partial_{x_1z} \tilde n$. This operator is elliptic,
since the associated matrix has only positive eigenvalues. Looking
after travelling wave solutions, we search a speed $c$ and a
profile $u(x,z)$ such that  $\tilde n(t,x,z):=u(x\cdot
e_1-ct,z)=u(x_1-ct,z)$ solves \eqref{eq1}. For convenience, we
drop the numerical subscript and write $x$ instead of $x_1$.
Hence, we are looking after $(c,u(x,z))$ such that
\begin{equation}
-\mathcal E(u)(x,z)-c u_{x}(x,z)= \left(r(z)- \displaystyle\int
_\R k(z,z')u(x,z')\,dz'\right)u(x,z) \quad \text{ in } \R
^2,\label{eq-tw}
\end{equation}
where
$$
\mathcal E (u):=u_{x x}+(B^2+1)u_{zz}-2B u_{xz}.
$$

\vskip 5pt Throughout the paper, we make the following assumption.
\begin{ass}[Structure of $r$ and $k$]\label{hyp} Function $r$ is in the H\"older space $C^{0,\theta}_{loc}(\R)$ for some $0<\theta<1$, and there is $\delta>0$ such
that
\begin{equation}\label{hyp-r}
\forall z\in \mathbb R, \quad r(z)\leq \frac 1\delta-\delta z^2.
\end{equation}
Function $k$ is in the H\"older space $C^{0,\theta}_{loc}(\R^2)$
and there are $k^->0$, $k^+>0$ such that
$$\forall (z,z')\in \R ^2,\quad  k^-\leq k(z,z')\leq  k^+.$$
\end{ass}

Let us next introduce a principal eigenvalue problem that is
necessary to enunciate our main result. For more details on
principal eigenvalue problems in general domains we refer to
\cite{ber-nir-var}, \cite{Ber-Ros} and the references therein.

\begin{defi}[Principal eigenvalue problem]\label{def:vp} We denote by
$\left(\lambda_\infty^0,\Gamma_\infty^0\right) \in \R \times
C^\infty(\R)$ the solution of the
 principal eigenvalue problem
\begin{equation}\label{vp-pb}
\begin{cases} -\left(B^2+1\right) \Delta _z\Gamma_\infty^0(z) -r(z)\Gamma_\infty^0(z) =\lambda_\infty^0 \Gamma_\infty^0(z)  \quad\text{ for all  } z\in \R\\
\Gamma_\infty^0(z) >0\quad \text{ for all } z\in \R,\quad
\Gamma_\infty^0(0)=1.
 \end{cases}
 \end{equation}
\end{defi}

Observe that in the case where $r(z)=1-Az^2$, $A>0$, we have
$\lambda_\infty^0=\sqrt{A\left(B^2+1\right)}-1$ and
$\Gamma_\infty^0
(z)=\exp\left(-\sqrt{\frac{A}{B^2+1}}\frac{z^2}2\right)$ is a
Gaussian profile.
\bigskip

We first state that as soon as $\lambda _\infty ^0>0$, extinction
of the population occurs.

\begin{prop}[Extinction]\label{extinction}
Assume $\lambda_\infty^0>0$. For any initial population $n^0$ such
that
$$\left\|\frac{n^0(x,y)}{\Gamma_\infty^0(y-Bx)}\right\|_{L^\infty(\mathbb
R^2)}<\infty,$$ the solution of \eqref{eq0} with initial condition
$n^0$ goes extinct exponentially fast as $t\to\infty$:
\begin{equation*}
\left\|\frac{n(t,x,y)}{\Gamma_\infty^0(y-Bx)}\right\|_{L^\infty(\mathbb
R^2)}=O(e^{-\lambda_\infty^0t}).
\end{equation*}
\end{prop}

Next, we state our main result: as soon as $\lambda_\infty^0 <0$,
invasion waves exist. Precisely, the following holds.

\begin{theo}[Travelling waves]\label{th:tw} Assume $\lambda_\infty^0 <0$
and define
\begin{equation}\label{vitesse-min}
c^*:=2 \sqrt{\frac{-\lambda_\infty^0}{B^2+1}}.
\end{equation}
Then the following holds. \begin{description} \item $(i)$ For all
$c\geq c^*$, there exists a positive $u\in C^2(\R ^2)$ solution of
\begin{equation}\label{eq-dans-espace}
-\mathcal E(u)(x,z)-c u_x(x,z)= \left(r(z)- \displaystyle\int _\R
k(z,z')u(x,z')\,dz'\right)u(x,z) \quad \text{ in } \R ^2,
\end{equation}
with
\begin{equation}\label{gauche}
\nu \mathbf{1}
_{\left\{(x,z)\in(-\infty,0]\times[-\nu,\nu]\right\}} \leq u(x,z)
\leq C e^{-K z^2},
\end{equation}
for some $\nu>0$, $C>0$, $K>0$, and
\begin{equation}\label{cond-bords}
\|u(x,\cdot)\|_\infty\to_{x\to +\infty} 0,\quad \int _{\R}
u(x,z)\,dz\to_{x\to +\infty} 0.
\end{equation}
Additionally, when $c>c^*$, there exists $\mu<0$ such that
\begin{equation}\label{controle-exp}
u(x,z)\leq e^{\mu \left(\sqrt{B^2+1}\,x+\frac
B{\sqrt{B^2+1}}z\right)}\Gamma_\infty^0 (z).
\end{equation}

\item $(ii)$ When $0\leq c<c^*$, there is no positive solution of
\eqref{eq-dans-espace} such that $\liminf _{x\to +\infty}u(x,0)=0$
and $u(x,z)\leq \psi(z)$ for some $\psi \in L ^1(\R)$.
\end{description}
\end{theo}

\subsection{Comments}

{\bf On the extinction case.} The proof of Proposition
\ref{extinction} is elementary and we now give the proof. The
result is a consequence of the parabolic comparison principle
satisfied by the local equation
\begin{equation}
\partial_t\phi(t,x,y)-\Delta_{x,y}\phi(t,x,y) =r(y-Bx.e)\phi(t,x,y).\label{maj1}
\end{equation}
Indeed, one can check that $M
e^{-\lambda_\infty^0t}\Gamma_\infty^0(y-Bx.e)$ and $n(t,x,y)$ are
respectively a super- and a sub-solution of \eqref{maj1} with
ordered initial data (for $M$ large enough).

\vskip 5pt

\noindent{\bf On the construction of waves.} Let us first comment
on a major difficulty in the construction of travelling fronts.
When the competition term is replaced by a local (in $x$ and $y$)
density regulation, many techniques based on the comparison
principle --- such as some monotone iterative schemes
 or the sliding method \cite{Ber-Nir3}--- can be used to get {\it a priori} bounds, existence and monotonicity properties of  the solution.
Since integro-differential equations with a nonlocal competition
term  do not satisfy the comparison
 principle, it is unlikely that such techniques apply here.

It turns out that the  considered problem here has some
similarities with the case of a population structured by a
 spatial variable only, that is $n=n(t,x)$, submitted to a nonlocal competition as studied in  \cite{Ber-Nad-Per-Ryz} (see also \cite{Alf-Cov}).
 In this work, the construction of a travelling front is  based on a sequence of approximating problems on intervals $(-a_n,a_n)$,
 with $a_n\to\infty$.  Due to the lack of comparison principle for the approximated problem, the construction of a solution is based on
  a  topological degree argument, a method introduced initially in \cite{Ber-Nic-Sch}.

To construct our fronts, we adopt a similar strategy  and consider
a sequence of   problems in growing  boxes
$(-a_n,a_n)\times(-b_n,b_n)$, with a normalization at the origin.
In order to make this strategy possible, a key point is to
establish \textit {a priori} estimates, independent on the size of
the boxes, on the profile $u$, the speed $c$ and, in particular,
the tails of $u$ when $z$ is large. Due to the nature of the
considered kernels here, an uniform estimate on $u$ is obtained
using a local pointwise $L^p$ estimate, whereas the uniform
control on $c$ is obtained by showing that our problem does not
have a solution if the speed $c$ is too large or if $c=0$. Notice
that the latter analysis may also be used to simplify the proof in
\cite{Ber-Nad-Per-Ryz}.

Let us highlight that, in contrast with \cite{Ber-Nic-Sch} and
\cite{Ber-Nad-Per-Ryz}, it is far from obvious that the
constructed travelling fronts are monotone w.r.t. $x$ for $x>0$
large enough. Therefore, we shall need an extra work to catch the
behavior \eqref{cond-bords} as $x\to+\infty$.

Notice also that the comprehension of the behavior of the wave as
$x\to-\infty$ is quite involved. In the related case of the
nonlocal Fisher-KPP equation, the positive steady state $u\equiv
1$ may present, for some kernels, a Turing instability (see e.g.
\cite{Gen-Vol-Aug}, \cite{Ber-Nad-Per-Ryz}, \cite{Alf-Cov}). Such
a situation may also occur in our context.

Let us also mention that, although we construct fronts without
relying on any monotonic properties of the profiles, it is
suspected, as in the case  of the nonlocal Fisher-KPP equation
\cite{Ber-Nad-Per-Ryz}, that there exist monotone and non monotone
travelling fronts. The understanding of such issues is quite
challenging.

\vskip 5pt

\noindent{\bf Organization of the paper.} In Section
\ref{intro_bio} we briefly describe the biological context of
 \eqref{eq0} and give an interpretation  of our results.  Then, we prove Theorem \ref{th:tw} in Sections \ref{s:box}--\ref{s:faster}. In Section \ref{s:box}, we start by deriving
  some  \textit{a priori} bounds and then, using  a Leray-Schauder topological degree argument, we construct a solution in a bounded
box. In Section \ref{s:minimal-speed}, we let the box tend to $\R
^2$ and  obtain a wave, which turns out to be the one with minimal
speed $c=c^*$. We also show the non existence of waves with speed
$0\leq c<c^*$. Lastly, we construct faster waves $c>c^*$ in
Section \ref{s:faster}.

\section{Biological interpretation of the results}\label{intro_bio}

In this section we briefly precise the biological context of
\eqref{eq0}.

In the present paper, we are interested in biological invasions
involving darwinian evolution. Species invading new territories
often face environmental gradients of e.g. temperature,
luminosity, antibiotic chemicals. Experimentally, it is well
documented that invasive species then evolve during their range
expansion \cite{Etterson}, \cite{Keller}, to adapt to local
conditions. To understand the speed, or even the success of an
invasion, one should thus consider the dispersion, birth and death
processes, but should also take into account evolution
\cite{Griffith}, \cite{Keymer}, \cite{Hermsen}. Those questions
become especially important in the context of the {\it global
warming} \cite{Davis}, \cite{Duputie}: the favorable environmental
conditions of many species move towards the north, implying
important changes in species' range. It is also of great
importance for the evolution of resistance of bacteria to
antibiotics \cite{Hermsen}.

More generally, many evolutionary biology questions involve
spatially structured populations, while most existing models
 either neglect the spacial structure of the population, or largely simplify it. New theoretical tools are then needed,
 and structured population models are natural candidates: they enable the modelling of all the biological phenomena mentioned
  above, and numerical simulations show that they are able to reproduce interesting features. Analyzing this type of model is
  however challenging, even in a homogeneous setting, see e.g. \cite{Jabin}, \cite{Lorz}. This work, as well as the results of \cite{Bouin}
   are first steps in the mathematical understanding of the
dynamics of those models.

 The main application of our result concerns asexual populations living in an environmental cline. The simplest model then writes
\begin{eqnarray*}
&\partial_tn(t,x,y)&-\frac{\sigma_x^2}2\Delta_xn(t,x,y)-\frac{\sigma_m^2}2\Delta_{y}n(t,x,y)=\\
&&\left(r_{max}-\frac 1{2V_s}(y-bx)^2-\frac 1K\int
_{\R}n(t,x,y')\,dy'\right)n(t,x,y),
\end{eqnarray*}
 where $\sigma_x$, $\sigma_m$ describe respectively the diffusion rate and the mutation rate of the population, $\frac 1{2V_s}$ is
  the strength of the selection, $b$ is the steepness of the environment cline, and $K$ the carrying capacity of the environment.
After the rescaling $n(t,x,y)=\tilde
n(r_{max}t,\frac{\sqrt{2r_{max}}}{\sigma_x}x,\frac{\sqrt{2r_{max}}}{\sigma
_m}y)$, we see that $\tilde n$ solves
 \eqref{eq0} with
 $$
 r(y)=1-Ay^2,\quad A=\frac{\sigma_m^2}{4r_{max}^2 V_s},\quad B=\frac{\sigma_x }{\sigma_m}b,\quad k\equiv \frac 1{K r_{max}}.
 $$
The population then gets extinct if $A\left(B^2+1\right)>1$, while
if $A\left(B^2+1\right)<1$, invasion fronts exist, with a
 minimal propagation speed (in the original variables)
  $$
  \sqrt{2 r_{max}}\sigma_x\left(1-\frac{\sigma_m}{2r_{max}\sqrt{V_s}}\sqrt{\left(b\frac{\sigma_x}
  {\sigma_m}\right)^2+1}\right)^{1/2}\left(\left(b\frac{\sigma_x}{\sigma_m}\right)^2+1\right)^{-1/2}.$$

\begin{rem}
There exists thus only two dynamics: either the population gets
extinct, or it succeeds to invade the whole territory. The
situation of asexual populations is then very different from the
case of sexual populations (see \cite{Kirkpatrick},
\cite{Mir-Rao}), where populations surviving with a limited range
only are possible.

 Notice also that during invasions, the dispersion of individuals can evolve (see e.g. \cite{Phillips}). Our result does not
  apply to this problem, and we refer to \cite{Benichou}, \cite{Bouin} for such a
  situation.
  \end{rem}

\section{The problem in a bounded box}\label{s:box}

\subsection{On some principal eigenvalue problems}\label{ss:eigenvalues}

We first introduce some principal eigenvalue problems, whose
eigenfunctions will serve as boundary conditions when stating the
travelling wave problem in a bounded box.

For $\nu\in [0,\delta)$, where $\delta>0$ is as in Assumption
\ref{hyp}, we denote by $\left(\lambda
_\infty^{\nu},\Gamma_\infty^{\nu}\right)$ the solution of the
 principal eigenvalue problem
\begin{equation}\label{vp-delta-pb}
\begin{cases} -\left(B^2+1\right) \Delta _z\Gamma_\infty^{\nu}(z) -\left(r(z)+\nu z^2\right)\Gamma_\infty^{\nu}(z) =\lambda _\infty^{\nu} \Gamma_\infty^{\nu}(z)  \quad\text{ for all  } z\in \R\\
\Gamma_\infty^{\nu}(z) >0\quad \text{ for all } z\in \R,\quad
\Gamma_\infty^{\nu}(0)=1.
 \end{cases}
 \end{equation}
Notice that this definition is coherent with \eqref{vp-pb}, and
for any $\nu\in [0,\delta)$, we have $\lambda _\infty^{\nu} \leq
\lambda_\infty^0 <0$. Also, for $\nu\in [0,\delta)$ and $b>0$, we
define $\left(\lambda_b^{\nu},\Gamma_b^{\nu}\right)$ as the
solution of the  principal eigenvalue problem
\begin{equation}\label{vp-delta-pb-moins}
\begin{cases} -\left(B^2+1\right)\Delta_z\Gamma _b^{\nu}(z)-\left(r(z)+\nu z^2\right)\Gamma _b^{\nu}(z)=\lambda _b^{\nu} \Gamma _b^{\nu}(z) \quad\text{ for all } z\in(-b,b)\\
\Gamma_b^{\nu}(\pm b)=0\\
 \Gamma_b^{\nu}(z) >0\quad \text{ for all } z\in
 (-b,b),\quad \Gamma_b^{\nu}(0)=1.
 \end{cases}
 \end{equation}
Let us observe that  $b<b'$ implies
$\lambda_\infty^{\nu}<\lambda_{b'}^{\nu}<\lambda_b^{\nu}$, and
that $\lambda_b^{\nu}\to \lambda_\infty^{\nu}$ as $b\to\infty$. To
construct the travelling waves, we will use the eigenfunctions
$\Gamma_b^{\delta/3}$, for $b>0$ as a boundary value. To bound
from above those functions independently of $b>0$, we will also
use the functions $\Gamma^{2\delta/3}_\infty$. Notice that
$-\max_{z\in\mathbb R} r(z)\leq \lambda^{2\delta/3}_\infty <
\lambda_\infty^{\delta/3}<0$. To show that
$\Gamma^{2\delta/3}_\infty$ is integrable, we define $C:=\max
_{[-\bar z,\bar z]} \Gamma^{2\delta/3}_\infty$, where $\bar
z:=\frac {\sqrt 6}{\delta}$, and
$\psi(z):=C\textrm{exp}\left(-\sqrt{\frac{\delta}{B^2+1}}\frac{z^2-\bar
z^2}{2\sqrt 6}\right)$, so that $\psi(\pm \bar z)=C$ and
$$
-\left(B^2+1\right)\Delta _z \psi(z)-\left(r(z)+\frac{2\delta}3
z^2\right)\psi(z)\geq \left(\frac \delta 6 z^2-\frac 1 \delta
\right)\psi(z)\geq 0,
$$
for all $z\in(-\infty,-\bar z)\cup(\bar z,\infty)$. Since in
$(-\infty,-\bar z)\cup(\bar z,\infty)$ we have
$r(z)+\frac{2\delta}3z^2\leq \frac 1\delta-\frac{\delta} 3z^2\leq
0$, the comparison principle then applies to \eqref{vp-delta-pb}
on $(-\infty,-\bar z)\cup(\bar z,\infty)$, and yields
$\Gamma^{2\delta/3}_\infty (z)\leq \psi (z)$ for all
$z\in(-\infty,-\bar z)\cup(\bar z,\infty)$. As a result, for some
constant which we denote again by $C$, we have
\begin{equation}\label{Gamma-bar-gaussian} \Gamma^{2\delta/3}_\infty (z) \leq
C\textrm{exp}\left(-\sqrt{\frac{\delta}{B^2+1}}\frac{z^2}{2\sqrt
6}\right) \quad\text{  for all } z\in \R,
\end{equation}
which implies in turn that $ \Gamma^{2\delta/3}_\infty \in
L^1(\R)$.

For a given $b_0>0$, we now use a similar argument to control the
functions $\Gamma _b^{\delta/3}(z)$ uniformly w.r.t.
$b\in[b_0,\infty]$. When $z$ lies in $[-\bar z,\bar z]$, where
$\bar z:=\max\left\{\sqrt{\frac{3\left(\lambda _{b_0}^{\delta/3} -
\lambda_\infty^{2\delta/3}\right)}\delta},\sqrt{\frac
3{2\delta}\left(\lambda _{b_0}^{\delta/3}+\frac 1
\delta\right)}\right\}$, the coefficients of the equations in
\eqref{vp-delta-pb} and \eqref{vp-delta-pb-moins} are uniformly
bounded w.r.t. $b\in[b_0,\infty]$. Therefore the Harnack
inequality implies that there is $C>0$ such that $\Gamma
_b^{\delta/3}(z)\leq  C$, for all $z\in[-\bar z,\bar z]$, all
$b\in[b_0,\infty]$. By the definition of $\bar z$ we see that, on
the one hand, $ \Gamma^{2\delta/3}_\infty$ is a super-solution for
\eqref{vp-delta-pb} and \eqref{vp-delta-pb-moins} --- with
$\nu=\frac\delta 3$--- on $(-\infty,-\bar z)\cup(\bar z,\infty)$
and that, on the other hand, the comparison principle applies.
Therefore, there exists $\bar C>0$ such that
\begin{equation}\label{controle-Gamma-b}
\Gamma_b^{\delta/3}(z)\leq \bar C\Gamma^{2\delta/3}_\infty(z)\leq
\bar C  \Vert \Gamma^{2\delta/3}_\infty \Vert _\infty \quad\text{
for all } z\in \R,\, b\in[b_0,\infty].
\end{equation}
In particular, we have
\begin{equation}\label{controle-mass-bords}
\int _{\R} \Gamma_b^{\delta/3}(z)\,dz\leq \mass :=\bar C\int _{\R}
\Gamma^{2\delta/3}_\infty(z)\,dz<\infty\quad\text{ for all }
b\in[b_0,\infty].
\end{equation}

\subsection{The problem in a box}\label{ss:box}

 For $a>0$, $b>0$ and $\ep
\in(0,1)$, we consider the problem of finding a speed $c\in \R$
and a real function $u(x,z)$, defined for $(x,z)\in
[-a,a]\times[-b,b]$, such that
$$
P(a,b,\ep)\quad\begin{cases}\, -\mathcal E(u)(x,z)-c u_x(x,z)\vspace{3pt}\\
\quad=\mathbf{1} _{\{u(x,z)\geq 0\}} \left(r(z)- \displaystyle\int
_{-b}^bk(z,z')u(x,z')\,dz'\right)u(x,z) \quad &\text{ in }Q\vspace{8pt} \\
\,u(x,z)=\mathbf{1}_{\{x=-a\}}(x)\Gamma_b^{\delta/3}(z)\quad &\text{ on }\partial Q\vspace{5pt}\\
\, u(0,0)=\ep,
\end{cases}
$$
where $Q:=(-a,a)\times(-b,b)$.  The elliptic operator is given by
\begin{equation}\label{def-elliptique}
-\mathcal E (u)=-u_{xx}-\left(B^2+1\right)u_{zz}+2Bu_{xz}.
\end{equation}

If $(c,u)$ is a solution achieving a negative minimum at
$(x_m,z_m)$ then, from the boundary conditions we deduce that
$(x_m,z_m)$ lies in the interior of the rectangle, and that
$-\mathcal E(u)-c u_x=0$ on a neighborhood of $(x_m,z_m)$. The
maximum principle thus implies $u\equiv u(x_m,z_m)$, which cannot
be. Therefore any solution of $P(a,b,\ep)$ satisfies $u\geq 0$
and, by the strong maximum principle,
\begin{equation}
u> 0\,\quad\text{and }\quad   -\mathcal E(u)(x,z)-c u_x(x,z)=\
\left(r(z)- \int _{-b}^bk(z,z')u(x,z')\,dz'\right)u(x,z) \quad
\text{ in } Q\,. \label{alternative}
\end{equation}

\medskip
In the following of the section, we shall construct a solution to
$P(a,b,\ep)$ via a Leray-Schauder topological degree argument. To
make this possible, we consider a family of problems as follows.
For $a>0$, $b>0$ and $\tau\in[0,1]$, we consider the problem of
finding a speed $c\in \R$ and a nonnegative real function $u(x,z)$
such that
$$
P_\tau(a,b)\begin{cases}\, -\mathcal E(u)(x,z)-c u_x(x,z)\vspace{3pt}\\
\quad=\left(r(z)- \tau\displaystyle\int
_{-b}^bk(z,z')u(x,z')\,dz'-\gamma(1-\tau)u(x,z)\right)u(x,z) \quad &\text{ in }Q\vspace{8pt} \\
\,u(x,z)=\mathbf{1}_{\{x=-a\}}(x)\Gamma_b^{\delta/3}(z)\quad
&\text{ on }\partial Q,
\end{cases}
$$
where $\gamma>0$ will be specified later (see Lemma
\ref{lem:local-pb}). Note that $P_0(a,b)$ reduces to a local
problem, and that solving $P(a,b,\ep)$ is equivalent to solving
$P_1(a,b)$ with the additional normalization condition $
u(0,0)=\ep$.

\begin{rem} A first natural idea to define a family of problems
would be to consider $-\mathcal E(u)-c u_x=\tau \left(r(z)- \int
_{-b}^bk(z,z')u(x,z')\,dz'\right)u$. But then we cannot get a
uniform w.r.t. $0\leq \tau \leq 1$ control of the tails of $u$
(see Lemma \ref{lem:apriori-tails}), which is crucial to derive
e.g. a lower bound on the standing waves (see Lemma
\ref{lem:a-priori-speed2}). This is the reason why we consider the
family $P_\tau(a,b)$ as above. Therefore, the topological degree
argument (see subsection \ref{ss:construction-box}) is rather
involved and requires to analyze a whole family of local problems
(see Lemma \ref{lem:locsigma}).
\end{rem}

\subsection{A priori estimates for $u$}\label{ss:apriori-u}

We provide a priori bounds for the profile $u$ of solutions to
$P_\tau(a,b)$. When $0\leq \tau \leq 1/2$, the local part of the
equation shall be enough to derive Lemma \ref{lem:apriori-bound}.
On the other hand, when $1/2\leq\tau\leq 1$, the nonlocal part is
quite relevant and we first need to control the vertical mass of
$u$, namely
$$
m(x):=\int _{-b} ^b u(x,z)\,dz.
$$

\begin{lem}[A priori bound for the mass]\label{lem:apriori-mass}  For all $a>0$, $b\geq b_0>0$, $1/2\leq \tau \leq 1$, any
 solution $(c,u)$ of $P_\tau(a,b)$ satisfies
 $$
  0 \leq \int _{-b}^b u(x,z)\,dz \leq \max\left(\frac{2\max_{\mathbb R}r}{k^-},\mass\right),\quad \forall
  x\in[-a,a].
  $$
\end{lem}

\noindent {\bf Proof.} Integrating w.r.t. $z$ the inequality
$-u_{xx}-cu_x\leq \left(\max_{\mathbb R}r-\frac 1 2 k^-
\,m(x)\right)u+\left(B^2+1\right)u_{zz}-2Bu_{xz}$, we get
\begin{eqnarray*}
 -m''(x)-cm'(x)&\leq& \left(\max_{\mathbb R}r-\frac 1 2 k^- \,m(x)\right)m(x)+\left(B^2+1\right)\left(u_z(x,b)-u_z(x,-b)\right)\\
&&-2B\left(u_x(x,b)-u_x(x,-b)\right).
\end{eqnarray*}
Since $u _z(x,b)\leq 0$, $u_z(x,-b)\geq 0$ and
$u_x(x,b)=u_x(x,-b)=0$, the mass satisfies the Fisher-KPP
inequality $-m''-cm'\leq \left(\max_{\mathbb R}r-\frac 1 2 k^-
\,m\right)m$. Since $m(-a)\leq \mass$ (see
\eqref{controle-mass-bords}) and $m(a)=0$, the maximum principle
concludes the proof of the lemma. \fdem

\medskip The above nonlocal control now provides the following a priori
bound for $u$.

\begin{lem}[A priori bound for $u$]\label{lem:apriori-bound} There exists $M>0$ such that, for all $a>0$, $b\geq
b_0>0$, $0\leq \tau \leq 1$, any solution $(c,u)$ of $P_\tau(a,b)$
with $0\leq c\leq c^\ast+1$ satisfies
$$0 \leq u(x,z) \leq M,\quad \forall
(x,z)\in[-a,a]\times[-b,b].$$
\end{lem}

\noindent {\bf Proof.}  For $0\leq \tau \leq 1/2$, one keeps the
local part and writes $-\mathcal E (u)-cu_x\leq \left(\max _\R r
-\frac \gamma 2 u\right)u$; recalling \eqref{controle-Gamma-b},
the maximum principle then gives a control of $u$ by
$\max\left(\frac {2 \max r} \gamma,\bar C \left\Vert
\Gamma^{2\delta/3}_\infty \right\Vert _\infty\right)$.

Next, for $1/2\leq \tau \leq 1$, let us denote by $(x_M,z_M)$ a
point where $u$ achieves its maximum $M$, and by $B_r$ the ball
centered at $(x_M,z_M)$ with radius $r>0$. Note that $0\leq u\leq
\bar C\left\|\Gamma^{2\delta/3}_\infty\right\|_\infty$ on
$\partial Q$. The function $w:=u-\bar
C\left\|\Gamma^{2\delta/3}_\infty\right\|_\infty$ therefore
satisfies
$$
\begin{cases}\, -\mathcal E(w)-c w_x-(\max _{\R} r)w\leq C_0:=(\max _\R r)\bar C\left\|\Gamma^{2\delta/3}_\infty\right\|_\infty \quad &\text{ in }Q\vspace{4pt} \\
\,w\leq 0 &\text{ on } \partial Q.
\end{cases}
$$
{}From the local maximum principle  \cite[Theorem 9.20]{Gil-Tru}
and it extension up to balls intersecting the boundary of the
domain \cite[Theorem 9.26]{Gil-Tru}, we infer that
$$\sup_{B_{1/2} \cap Q} w \le C_1\left( \frac{1}{|B_1|}\int_{B_1\cap Q}w^+ + C_2 \Vert C_0 \Vert _ {L^N(B_1\cap
Q)}\right),
$$
where $C_1=C_1(B)$ and $C_2=C_2(B)$ are positive constants. Notice
that $C_1$ does not depend on $c$, which is a coefficient of the
operator $L(w):=\mathcal E (w)+cw_x+w$, because $c$ belongs to a
bounded interval, namely $[0,c^\ast+1]$. Using successively $u\geq
0$ and Lemma \ref{lem:apriori-mass} we deduce that
$$\int_{B_1\cap Q}w^+\le \int_{B_1\cap Q}u +  \int_{B_1\cap Q}\bar C\left\|\Gamma^{2\delta/3}_\infty\right\|_\infty \le
 2\max\left(\frac{2\max_{\mathbb R}r}{k^-},\mass\right)+ \bar C\left\|\Gamma^{2\delta/3}_\infty\right\|_\infty|B_1|.
$$
Recalling that $M=\max u$ is achieved at the center of the ball
$B_{1/2}$, we deduce form the upper estimates that
$$M \le \bar C\left\|\Gamma^{2\delta/3}_\infty\right\|_\infty+ C_1\left( \frac{2}{|B_1|}\max\left(\frac{2\max_{\mathbb R}r}{k^-},\mass\right)
  +\bar C\left\|\Gamma^{2\delta/3}_\infty\right\|_\infty+ C_2 C_0 |B_1|\right).$$
This concludes the proof of the lemma. \fdem

\medskip

We now provide a control of the tails of the solutions as $|z|\to
\infty$ by appropriate Gaussian functions (recall estimate
\eqref{Gamma-bar-gaussian}).

\begin{lem}[Gaussian control of the tails of $u$]\label{lem:apriori-tails}
There exists $\bar M>0$ such that, for all $a>0$, $b\geq b_0>0$,
$0\leq \tau \leq 1$, any solution $(c,u)$ of $P_\tau(a,b)$ with
 $0\leq c\leq c^\ast+1$ satisfies
$$
0\leq u(x,z)\leq  \bar M \Gamma^{2\delta/3}_\infty(z),\quad
\forall (x,z)\in[-a,a]\times[-b,b].
$$

\end{lem}

\noindent {\bf Proof.} First observe that
\begin{equation}\label{eq:apriori-tails}
-\mathcal E(u)-c u_x-r(z)u\leq 0 \quad\text{on } Q,\quad u\leq
\bar C  \Gamma^{2\delta/3}_\infty \quad\text{on } \partial Q.
\end{equation}
Define $\bar\phi(x,z)=\bar \phi(z):=\bar M
\Gamma^{2\delta/3}_\infty(z)$, with $\bar M>0$ to be specified
later. Recall that $\Gamma^{2\delta/3}_\infty$ solves
\eqref{vp-delta-pb} so that $ -\e(\bar \phi) -c\bar\phi_x
-r(z)\bar \phi =\left(\frac{2\delta}3 z^2+
\lambda^{2\delta/3}_\infty\right)\bar
 \phi$. Therefore, if $\beta:=\max\left(\frac 1\delta,\sqrt{\frac{-3
\lambda^{2\delta/3}_\infty}{2\delta}}\right)$, we have
\begin{equation}\label{sur-sol-gaussienne}
-\e(\bar \phi) -c\bar\phi_x -r(z)\bar \phi\ge 0\quad\text{ on }
(-a,a)\times(\beta,b).
\end{equation}
Let us now select
\begin{equation}\label{choice}
\bar M:=\max\left( \bar C,\frac {M}
{\min_{[-\beta,\beta]}\Gamma^{2\delta/3}_\infty}\right),
\end{equation}
where $M$ is as in the previous lemma. The choice \eqref{choice}
enforces $\bar \phi (x,z)\geq \bar C \Gamma^{2\delta/3}_\infty(z)
\geq u(x,z)$ on $\{\mp a\}\times [\beta,b]\cup [-a,a]\times
\{b\}$, and $\bar \phi (x,z)\geq M\geq u(x,z)$ on
$[-a,a]\times\{\beta\}$. Hence the comparison principle
--- note that $r(z)\leq 0$, when $z\geq \beta$--- yields $u\leq
\bar \phi$ on $[-a,a]\times [\beta,b]$ and, by the choice
\eqref{choice}, on $[-a,a]\times[0,b]$. Similarly, $u\leq \bar
\phi$ on $[-a,a]\times [-b,0]$. The lemma is proved.
 \fdem

\subsection{A priori estimates for $c$}\label{ss:apriori-c}

We provide a priori bounds for the speed $c$ of solutions to
$P_\tau(a,b)$. We first show that, roughly speaking, too rapid
waves solutions of $P_\tau(a,b)$ have too small value at
$(x,z)=(0,0)$. We recall that the speed $c^*$ was defined in
\eqref{vitesse-min}.

\begin{lem}[A priori upper bound for $c$]\label{lem:a-priori-speed}
Let $b>0$ and $\ep \in(0,1)$ be arbitrary. Then there exists
$a_0=a_0(\ep,b)>0$ such that, for all $a\geq a_0$, all $0\leq \tau
\leq 1$, any solution $(c,u)$ of $P_\tau (a,b)$ with $c>c^*$
satisfies $u(0,0) <\ep$
--- and therefore cannot solve $P(a,b,\ep)$.
\end{lem}

\noindent {\bf Proof.} Let $b>0$ be given. Assume $c>c^*$ and let
us show that $u(0,0)\to 0$ as $a\to \infty$.

The function $u$  satisfies $-\mathcal E (u)-cu_x - r(z)u\leq 0$
in $Q$. Therefore, changing variables, the function
\begin{equation}\label{cgtvar}
 v(x,y):=u\left(\frac{x-By}{\sqrt{B^2+1}},\sqrt{B^2+1}\,y\right),
\end{equation}
satisfies
\begin{equation}\label{ineq-v}
 -v_{xx}-v_{yy}-c\sqrt{B^2+1}v_x-r\left(\sqrt{B^2+1}\,y\right)v\leq
 0,
\end{equation}
in $Q_1:=\left\{(x,y):\,|y|< \frac
b{\sqrt{B^2+1}},\,\left|\frac{x-By}{\sqrt{B^2+1}}\right|<a\right\}$.

We shall now construct a positive solution. Since $c>c^*$, one can
select $\mu<0$ such that $\mu ^2+c\sqrt{B^2+1}\mu+\frac{{c^*}^2}4
\left(B^2+1\right)=0$ and define $\varphi(s):=e^{\mu s}$ which
solves
$$
-\varphi ''-c\sqrt{B^2+1}\varphi '-\frac{{c^*}^2}4
 \left(B^2+1\right)\varphi =0.
 $$
Now,  let us define
$$
w (x,y):=\kappa_a \varphi(x)\Gamma_\infty^0\left(\sqrt{B^2+1}
y\right),\quad
\kappa_a:=\frac{\left\|\Gamma_b^{\delta/3}\right\|_\infty}{\min_{[-b,b]}\Gamma_\infty^0}
e^{\mu\left(a\sqrt{B^2+1}-\frac{Bb}{\sqrt{B^2+1}}\right)},
$$
with $\Gamma_\infty^0$ the eigenfunction appearing in Definition
\ref{def:vp}. Using direct computations and the definition of
$c^*$ in \eqref{vitesse-min} we see that
\begin{equation}\label{sol-positive}
-w_{xx}-w_{yy}-c\sqrt{B^2+1}w_x-r\left(\sqrt{B^2+1}\,y\right)w = 0
\quad\text{ in } Q_1.
\end{equation}

We now compare the values of $v$ and $w$ on the boundary. If
$(x,y)\in\partial Q_1$ is such that $\frac{x-By}{\sqrt{B^2+1}}\neq
-a$, then $v(x,y)=0< w(x,y)$. If $(x,y)\in\partial Q_1$ is such
that $\frac{x-By}{\sqrt{B^2+1}}= -a$ then
$$
v(x,y)=\Gamma_b^{\delta/3}\left(\sqrt{B^2+1}\,y\right)\leq
\kappa_{a}
e^{\mu\left(-a\sqrt{B^2+1}+\frac{Bb}{\sqrt{B^2+1}}\right)}\min_{[-b,b]}\Gamma_\infty^0
\leq w(x,y).
$$
As a result, we have $v\leq w$ on $\partial Q_1$.

Now a classical argument will imply $v\leq w$ on the whole of
$\overline{Q_1}$. Indeed, since $v$ is bounded and $w>0$ on
$\overline{Q_1}$, we can define
$\alpha_0:=\max\left\{\alpha>0:\,\alpha v\leq w\text{ in }
\overline{Q _1}\right\}>0$.  Then $\alpha _0v\leq w$ and there is
a point $(x_0,y_0)$ such that $\alpha _0v(x_0,y_0)=w(x_0,y_0)$. In
view of \eqref{ineq-v}, \eqref{sol-positive} and the strong
maximum principle the point $(x_0,y_0)$ has to lie on $\partial
Q_1$, which enforces $\alpha _0\geq 1$. Thus $v\leq w$ in
$\overline{Q_1}$, and
$$u(0,0)=v(0,0)\leq w(0,0)= \kappa_{a}\to 0\textrm{ as }a\to\infty,$$
which concludes the proof of the lemma.\fdem

\medskip

Next, we show that standing waves (i.e. $c=0$) have too large
value at $(x,z)=(0,0)$.

 \begin{lem}[Standing waves: a priori lower bound for $u(0,0)$]\label{lem:a-priori-speed2}
 There is $\ep ^*>0$ such that if $a,\,b$ are large enough, then, for all $0\leq \tau \leq 1$, any
standing solution $(c=0,u)$ of $P_\tau(a,b)$ satisfies $u(0,0)>\ep
^*$ --- and therefore cannot solve $P(a,b,\ep)$ for any $\ep\in
(0,\ep^\ast)$.
\end{lem}

\noindent {\bf Proof.} For $R>0$, let us introduce
$(\mu_R,\Upsilon_R)$ as the solution of the  principal eigenvalue
problem
\begin{equation}\label{vpboxR}
\begin{cases} -\mathcal E (\Upsilon_R)(x,z)-r(z)\Upsilon_R(x,z)=\mu_R \Upsilon_R(x,z)\quad \text{ for all } (x,z)\in(-R,R)^2\\
 \Upsilon_R = 0 \quad \text{ on } \partial \left((-R,R)^2\right)\\
  \Upsilon_R(x,z)>0\quad \text{ for all } (x,z)\in (-R,R)^2,\quad \Upsilon_R(0,0)=1.
 \end{cases}
 \end{equation}
If $R=\infty$, the above problem is equivalent to \eqref{vp-pb},
and $\mu_R\to\lambda_\infty^0$ as $R\to\infty$. Let us therefore
fix $R>0$ large enough so that
\begin{equation}\label{vpineq}
 \lambda_\infty^0 \leq \mu_R<\frac{\lambda_\infty^0}2<0,
\end{equation}
and
\begin{equation}\label{vpineq2}
k^+\int_{[-R,R]^c}\bar
M\Gamma^{2\delta/3}_\infty(z)\,dz\leq\frac{-\lambda_\infty^0}4.
\end{equation}

Next, let $a \geq R+1$, $b\geq R+1$, $0\leq \tau \leq 1$ be given,
and $(c=0,u)$ be a solution of $P_\tau(a,b)$. Thanks to
 the Harnack inequality, there exists $C>0$ (independent of $a$ and $b$), such that $\|u\|_{L^\infty\left([-R,R]^2\right)}\leq Cu(0,0)$, and then
$$
\max_{(x,z)\in [-R,R]^2}\left(\tau\displaystyle\int
_{-b}^bk(z,z')u(x,z')\,dz'+\gamma(1-\tau)u\right)\leq C
\left(2Rk^++\gamma\right)u(0,0)+\frac{-\lambda_\infty^0}4,
$$
where we have used Lemma \ref{lem:apriori-tails} and
\eqref{vpineq2}. As a result, $u$ satisfies
$$
-\mathcal E (u)-r(z)u\geq
\left(-C\left(2Rk^++\gamma\right)u(0,0)+\frac{\lambda_\infty^0}4\right)u
\quad\text{ in } (-R,R)^2,
$$
and $u>0$ on $\partial \left((-R,R)^2\right)$. Hence, if
$$
u(0,0)\leq
\frac{-\lambda_\infty^0/4}{C\left(2Rk^++\gamma\right)}=:\ep^*,
$$
$u$ becomes a super-solution for \eqref{vpboxR} --- thanks to
\eqref{vpineq}. We conclude as in the proof of Lemma
\ref{lem:a-priori-speed}: defining
$\alpha_0:=\max\left\{\alpha>0:\,\alpha \Upsilon_R \leq u \text{
in } [-R,R]^2\right\}>0$ and using the strong maximum principle we
see that $\alpha _0 \Upsilon_R \equiv u$ on $[-R,R]^2$, which
contradicts $u>0$. It follows that $u(0,0)>\ep^*$. The lemma is
proved. \fdem

\subsection{On some related local problems}\label{ss:local-pb}

First, we show the well-posedness of the local problem $P_0(a,b)$.

\begin{lem}[Well-posedness for $\tau =0$]\label{lem:local-pb} There exists $\gamma>0$ such that, if $a,\,b$ are large enough and $\ep\in (0,\ep^\ast)$,
there exists a unique $(c,u\geq 0)\in \R \times W^{2,\infty}(Q)$,
solution of $P_0(a,b)$, namely
$$
P_0(a,b)\begin{cases}\, -\mathcal E(u)(x,z)-c u_x(x,z)=
\left(r(z)-\gamma \,u(x,z)\right)u(x,z) \quad &\text{ in } Q\vspace{8pt} \\
\,u(x,z)=\mathbf{1}_{\{x=-a\}}(x) \Gamma_b^{\delta/3}(z)&\text{ on
}\partial Q,\end{cases}
$$
such that $u(0,0)=\ep$. Moreover, by the above a priori estimates,
 $0<c\leq c^* $ and $0\leq u \leq M$.
\end{lem}

\noindent {\bf Proof.} A standard argument proves that there is a
unique positive solution to $P_0(a,b)$. For the convenience of the
reader let us prove this fact. Since $0$ and a large enough
positive constant are respectively a sub- and a super-solution of
$P_0(a,b)$, the existence of a positive solution to $P_0(a,b)$ can
be obtained using a classical monotone iterative scheme. Also, by
the maximum principle, any positive solution of $P_0(a,b)$ is
bounded. Now let $u$ and $v$ be two bounded positive solutions of
$P_0(a,b)$. Thanks to the boundary condition and the Hopf lemma
the following quantity is well defined:
$$
\tau^*:=\inf \left\{\tau >0:\, \forall (x,z) \in \overline Q,
u(x,z)\le \tau v(x,z) \right\}.
$$
Assume  by contradiction that $\tau^*>1$. From the definition of
$\tau^*$, the boundary condition and the Hopf lemma, there exists
$(x_0,z_0)\in Q$ such that $u(x_0,z_0)=\tau^*v(x_0,z_0)$. At this
point, we get the contradiction
$$0\le \e(\tau^* v-u)(x_0,z_0)+c (\tau^*v-u)_x(x_0,z_0)\le \gamma\tau^*(1-\tau^*) v^2(x_0,z_0)<0. $$
Thus $\tau^*\le 1$ and we have $u\le v$. By interchanging the role
of $u$ and $v$,  we get $u\equiv v$. Hence there is a unique
positive solution to $P_0(a,b)$.

In order to apply a sliding method, we slightly modify the
Dirichlet boundary conditions and consider, for small $\eta >0$,
$$
P_0^\eta(a,b)\begin{cases}\, -\mathcal E(v^\eta)(x,z)-c
v^\eta_x(x,z)=
\left(r(z)-\gamma\,v^\eta(x,z)\right)v^\eta(x,z) \quad &\text{ in } Q\vspace{8pt} \\
\,v^\eta(x,z)=\displaystyle\frac{x-a}{-2a}\Gamma_{b}^{\delta/3\,,\eta}(z)
\quad &\text{ on } \partial Q,\end{cases}
$$
where
$\Gamma_{b}^{\delta/3\,,\eta}(z):=\Gamma_{b}^{\delta/3}(z)+\eta\Gamma_\infty^{\delta/3}(z)$.
Note that for any nonnegative $\eta$, $0$ and a large enough
positive constant are respectively a sub- and a super-solution of
$P_0^{\eta}(a,b)$. Thus the existence of a positive and bounded
solution of $P_0^{\eta}(a,b)$ can be obtained using a classical
monotone iterative scheme. We shall now select $\gamma
>0$ so that $(x,z)\mapsto \Gamma_b^{\delta/3,\,\eta}(z)$ become a
super-solution for $P_0^\eta(a,b)$. Using \eqref{vp-delta-pb},
\eqref{vp-delta-pb-moins}, $\lambda _b^{\delta/3} \geq \lambda
_\infty^{\delta/3}$ and $\Gamma _b^{\delta/3,\,\eta} \geq \Gamma
_b^{\delta/3}$ we see that
\begin{equation}\label{bords-sursol}
 -\left(B^2+1\right)\Delta_z\Gamma_b^{\delta/3,\,\eta}-\left(r(z)-\gamma\Gamma_b^{\delta/3,\,\eta}\right)\Gamma_b^{\delta/3,\,\eta}
\geq \left(\lambda_\infty^{\delta/3}+\frac \delta 3
z^2+\gamma\Gamma_b^{\delta/3}\right)\Gamma_b^{\delta/3,\,\eta},
\end{equation}
which is clearly nonnegative for $|z|\geq \bar
z:=\sqrt{\frac{-3\lambda_\infty^{\delta/3}}\delta}$. Now, as soon
as $b\geq \bar z +1$, it follows from the Harnack inequality that
there is $C>0$ such that $\Gamma_b^{\delta/3}(z)\geq \frac 1 C
\Gamma_b^{\delta/3}(0)=\frac 1 C$ for all
 $|z|\leq \bar z$. Hence if we select $\gamma:=-C\lambda _\infty^{\delta/3}>0$,
 the right-hand side member of \eqref{bords-sursol} becomes nonnegative
 for $|z|\leq \bar z$, and $(x,z)\mapsto
\Gamma_b^{\delta/3,\,\eta}(z)$ is a super-solution for
$P_0^\eta(a,b)$. Next, we show that, any solution of
$P_0^\eta(a,b)$ satisfies
\begin{equation}\label{ordre-strict}
0<v^\eta (x,z)<\Gamma_b^{\delta/3,\,\eta} (z),\quad \forall
(x,z)\in Q.
\end{equation}
Indeed, for $v^\eta$ a non negative solution  of
$P^{\eta}_0(a,b)$, the following quantity is well defined
$$\alpha _0:=\sup\{\alpha\geq 0:\, \forall (x,z)\in \bar Q,\, \alpha v^\eta(x,z)\le  \Gamma_b^{\delta/3,\,\eta} (z) \}\in(0,1],$$
since $\Gamma_b^{\delta/3,\,\eta}>0$ in  $\bar Q$. Let us assume
by contradiction that $\alpha _0<1$. In view of the boundary
conditions for $v^\eta$, this implies that a point $(x_0,z_0)$
where $\alpha _0 v^\eta(x_0,z_0)=\Gamma_b^{\delta/3,\,\eta} (z_0)$
cannot be on $\partial Q$. Hence
$w:=\Gamma_b^{\delta/3,\,\eta}-\alpha_0 v^\eta$ has a zero minimum
at $(x_0,z_0)\in Q$ and
$$
0\geq (-\e(w)-cw_x)(x_0,z_0)\geq \gamma \alpha _0
(v^\eta)^2(x_0,z_0)(1-\alpha _0)>0,
$$
which is absurd. Hence $\alpha _0=1$, and $0\le v^\eta\le
\Gamma_b^{\delta/3,\,\eta}$ in $\bar Q$.  Then
\eqref{ordre-strict} follows by applying  the strong maximum
principle.

By the classical sliding method \cite{Ber-Nir3}, $v^\eta$ is
strictly decreasing in the $x$ variable. For the convenience of
the reader, let us give a  proof of this fact. For $h>0$, define
$v^\eta_h(x,z):=v^\eta(x+h,z)$. Since $v^\eta_{2a}\equiv 0 <
v^\eta$ thanks to \eqref{ordre-strict}, one can define
$$
h^*:=\inf \left\{h>0:\, \forall \tau \in [h,2a], v^\eta_\tau \leq
v^\eta\right\}.
$$
Assume $h^*>0$. Then there are sequences $h_n \nearrow h^*$,
$(x_n,z_n)$ with $v^\eta_{h_n}(x_n,z_n)>v^\eta(x_n,z_n)$. After
extraction and using that the infimum $h^*$ is achieved, we have a
point $(x_\infty,z_\infty)$ such that
$v^\eta_{h^*}(x_\infty,z_\infty)=v^\eta(x_\infty,z_\infty)$, i.e.
a point of zero maximum for $v^\eta_{h^*}-v^\eta$. Because of the
boundary conditions the point $(x_\infty,z_\infty)$ cannot lie on
the upper or the lower boundary of $(-a,a-h^*)\times(-b,b)$. In
view of \eqref{ordre-strict}, it is neither allowed to lie on the
left or right boundary of $(-a,a-h^*)\times(-b,b)$. Since $v^\eta$
and $v^\eta_{h^*}$ are both solutions of $-\mathcal E(v)-c v_x=
\left(r(z) -\gamma v\right)v$ in $(-a,a-h^*)\times(-b,b)$, the
maximum principle then yields $v^\eta_{h^*}\equiv v^\eta$, i.e.
$v^\eta(x,z)=v^\eta(x+h^*,z)$. Hence
$v^\eta(x,z)=v^\eta(x+nh^*,z)$ for all $n\geq 0$. Letting $n\to
\infty$ yields $v^\eta\equiv 0$, which is a contradiction. It
follows that $h^*=0$ and $v^\eta$ is non increasing in the $x$
variable.

Now, we construct a solution to $P_0(a,b)$ as a limit, as $\eta
\to 0$, of solutions to $P_0^\eta(a,b)$. The interior elliptic
estimates imply that, for all $1<p<\infty$, the sequence
$(v^\eta)$ is bounded in $W^{2,p}(Q)$. From Sobolev embedding
theorem, one can extract a subsequence $(v^\eta)$ converging to
some $u$, strongly in $C^{1,\beta}(Q)$ and weakly in $W^{2,p}(Q)$.
Moreover, $u$ is a solution of $P_0(a,b)$. As a limit of
decreasing functions, $u$ is decreasing in the $x$ variable. By
differentiating the equation and applying the maximum principle,
one then obtains the strict decreasing of $u$ w.r.t. $x$.

It is then standard that, if $(c_1,u_1)$ and $(c_2,u_2)$ are two
solutions of $P_0(a,b)$ with $c_1>c_2$, then $u_1<u_2$. Indeed,
$u_2$ is a super-solution of the equation for $(c_1,u_1)$. Hence
there exists a solution for this equation which is below $u_2$. By
uniqueness  this solution is $u_1$. Hence $u_1\leq u_2$ and, by
the strong maximum principle, $u_1<u_2$. As seen in Lemma
\ref{lem:a-priori-speed}, Lemma \ref{lem:a-priori-speed2}, the
solution of $P_0(a,b)$ with speed $c=0$, $c>c^*$ satisfy $u(0,0)
>\ep^*>\ep$, $u(0,0)<\ep$ respectively, if $a,\,b$ are large
enough. Then, there is a unique $c$, which belongs to $(0,c^*]$,
such that the solution $(c,u)$ of $P_0(a,b)$ is $\ep$-normalized.
The lemma is proved. \fdem

\medskip In order to apply  a Leray-Schauder degree argument in
the next subsection, we also need to consider the family $0\leq
\sigma \leq 1$ of local problems
\begin{equation}\label{pblocsigma}
\tilde P_\sigma (a,b)\begin{cases}\, -\mathcal E(u)(x,z)-c
u_x(x,z)= \left(r(z)-(1-\sigma)\mathcal R
-\sigma\gamma u(x,z)\right)u(x,z) \quad &\text{ in }Q\vspace{8pt} \\
\,u(x,z)=\mathbf{1}_{\{x=-a\}}(x)\Gamma_b^{\delta/3}(z)\quad
&\text{ on }\partial Q,
\end{cases}
\end{equation}
where $\mathcal R:=\max_{z\in \R}\left(r(z)+\frac \delta 3
z^2\right)$.

\begin{lem}[On local problems $\tilde P_\sigma(a,b)$]\label{lem:locsigma}
\begin{description}
\item $(i)$ There exists $M>0$ such that, for all $a>0$, $b>0$,
$0\leq \sigma \leq 1$, any solution $(c,u\geq 0)$ of $\tilde
P_\sigma(a,b)$ satisfies $0\leq u \leq M$.

\item $(ii)$  Let $b>0$ and $\ep \in(0,1)$ be arbitrary. Then
there exists $a_0=a_0(\ep,b)>0$ such that, for all $a\geq a_0$,
all $0\leq \sigma \leq 1$, any solution $(c,u\geq 0 )$ of $\tilde
P_\sigma (a,b)$ with $c>c^*$ satisfies $u(0,0) <\ep$.

\item $(iii)$ There exists $\ep _0 >0$, $a_0>0$ such that, for any
$a\geq a_0$, there exists a speed $-\bar c=-\bar c(a)<0$ such that
for all $b\geq 1$, all $0\leq \sigma \leq 1$, any solution
$(c,u\geq 0)$ of $\tilde P_\sigma (a,b)$ with $c\leq -\bar c$
satisfies $u(0,0)
>\ep_0$.

\item $(iv)$ If $a,\,b$ are large enough and $\ep\in (0,\ep_0)$,
then there exists a unique $(c,u \geq 0)\in \R \times
W^{2,\infty}(Q)$ solution of
\begin{equation}\label{pblocsigmazero}
\tilde P_0 (a,b)\begin{cases}\, -\mathcal E(u)(x,z)-c u_x(x,z)=(r(z)-\mathcal R)u(x,z) \quad &\text{ in }Q\vspace{8pt} \\
\,u(x,z)=\mathbf{1}_{\{x=-a\}}(x)\Gamma_b^{\delta/3}(z)\quad
&\text{ on }\partial Q,
\end{cases}
\end{equation}
with $u(0,0)=\ep$. Moreover, by the above a priori estimates,
 $-\bar c<c\leq c^* $ and $0\leq u \leq M$.
\end{description}
\end{lem}

\noindent {\bf Proof.} Item $(i)$ follows from $-\mathcal E
(u)-cu_x\leq \sigma (\mathcal R  - \gamma  u)u$,
\eqref{controle-Gamma-b} and the maximum principle.

Define
$$
v(x,y):=u\left(\frac{x-By}{\sqrt{B^2+1}},\sqrt{B^2+1}\,y\right),\quad
Q_1:=\left\{(x,y):\,|y|< \frac
b{\sqrt{B^2+1}},\,\left|\frac{x-By}{\sqrt{B^2+1}}\right|<a\right\}.
$$
Then $v$ is a subsolution of \eqref{ineq-v} and, to prove $(ii)$,
we can reproduce the proof of Lemma \ref{lem:a-priori-speed}.

Let us prove $(iii)$. Observe that $v$ solves
\begin{equation}\label{eq-pour-v}
 -v_{xx}-v_{yy}-c\sqrt{B^2+1}v_x=  \left(r(\sqrt{B^2+1}\,y)-(1-\sigma)\mathcal R-\sigma \gamma
 v\right)v\quad \text{ in } Q_1,
\end{equation}
that $v(x,y)=0$ for $(x,y)\in\partial Q_1$ such that
$\frac{x-By}{\sqrt{B^2+1}}\neq-a$,
 and that $v(x,y)=\Gamma_b^{\delta/3}\left(\sqrt{B^2+1}\,y\right)$ for $(x,y)\in\partial Q_1$ such
that $\frac{x-By}{\sqrt{B^2+1}}= -a$. It follows from
\eqref{vp-delta-pb-moins} and the Harnack inequality that there
exists $C>0$ such that
\begin{equation}\label{harnack2}
\Gamma _b^{\delta/3}(z) >C \Gamma _b^{\delta/3}(0)=C,\quad \text{
for all } b\geq 1, \text{ all } |z|\leq 1.
\end{equation}
Define, for $\alpha > 0$,
$$
\psi_\alpha (x,y):=\frac C{\max_{[-1,1]} \Gamma^0_1
}\left(1-\frac{x+\alpha+\left(a\sqrt{B^2+1}+\frac{B}{\sqrt{B^2+1}}\right)}{2a\sqrt{B^2+1}}\right)
\Gamma^0_1\left(\sqrt{B^2+1}y\right),
$$
in $Q_2:=\left\{(x,y):\,|y|< \frac
1{\sqrt{B^2+1}},\,\left|\frac{x-By}{\sqrt{B^2+1}}\right|<a\right\}=Q_1
\cap \left\{|y|< \frac 1{\sqrt{B^2+1}}\right\}$. Since
$\psi_{\alpha=2a\sqrt{B^2+1}}\leq 0<v$ in $Q_2$, we can define
$$
\alpha_0:=\min\left\{\alpha\in\left[0,2a\sqrt{B^2+1}\right]:\,\forall
(x,y) \in \overline{Q_2}\,,\,\psi_{\alpha}(x,y)\leq
v(x,y)\right\}.
$$
Assume by contradiction that $\alpha _0>0$. Observe that
$\psi_{\alpha _0}(x,y)<0=v(x,y)$ for $(x,y)\in \partial Q_2$ such
that $\frac{x-By}{\sqrt{B^2+1}}=a$ and $|y|<\frac
1{\sqrt{B^2+1}}$, that $\psi_{\alpha _0}(x,y)=0<v(x,y)$ for
$(x,y)\in\partial Q_2$ such that $|y|=\frac 1{\sqrt{B^2+1}}$ and
$-a\leq \frac{x-By}{\sqrt{B^2+1}}<a$, and that --- by
\eqref{harnack2}--- $\psi_{\alpha
_0}(x,y)<\Gamma_b^{\delta/3}\left(\sqrt{B^2+1}\,y\right)=v(x,y)$
for $(x,y)\in\partial Q_2$ such that $\frac{x-By}{\sqrt{B^2+1}}=
-a$. Therefore the only two points on $\partial Q_2$ where
$v-\psi_{\alpha _0}$ may attain its zero minimum value are
$\left(\frac{a}{\sqrt{B^2+1}}\pm \frac B{\sqrt{B^2+1}},\pm\frac
1{\sqrt{B^2+1}}\right)$. But since $\alpha _0>0$ we see that, for
$\ep>0$ small enough, $\psi _{\alpha _0 -\ep}\leq 0\leq v$ in a
neighborhood of these two points of the boundary $\partial Q_2$.
Hence, by the definition of $\alpha _0$, there must be a point
$(x_0,y_0)\in Q_2$ where $v-\psi_{\alpha _0}$ attains its zero
minimum value, so that $0\geq -\Delta_{x,y}
\left(v-\psi_{\alpha_0}\right)(x_0,y_0)-c\sqrt{B^2+1}
\partial _x \left(v-\psi_{\alpha_0}\right)(x_0,y_0)$. Using
\eqref{eq-pour-v}, \eqref{vp-delta-pb} and straightforward
computations, we arrive at
\begin{eqnarray*}
0&\geq& v(x_0,y_0)\left(-(1-\sigma)\mathcal R-\lambda^0_1 -\sigma
\gamma v(x_0,y_0)\right)-\frac c{2a} \frac C {\max
 \Gamma^0_1}  \Gamma^0_1 \left(\sqrt{B^2+1}\,y_0\right)\\
&\geq & v(x_0,y_0)\left(-\mathcal R -\lambda^0_1 - \gamma
M\right)-\frac c{2a} \frac C {\max \Gamma^0_1}  \Gamma^0_1
\left(\sqrt{B^2+1}\,y_0\right).
\end{eqnarray*}
Since $v(x_0,y_0)=\psi_{\alpha _0}(x_0,y_0)\leq  \frac C {\max
\Gamma^0_1} \Gamma^0_1 \left(\sqrt{B^2+1}\,y_0\right)$ and since
$-\mathcal R -\lambda^0_1 \leq -\max_{[-1,1]}r- \lambda^0_1\leq 0$
we end up with
$$
0\geq \frac C {\max \Gamma^0_1} \Gamma^0_1\left(\sqrt{B^2+1}\,
y_0\right)\left(-\mathcal R -\lambda^0_1 -\gamma M -\frac
c{2a}\right).
$$
Hence, for large negative speed, namely
$$
c\leq -\bar c=-\bar c (a):=2a\left(-\mathcal R - \lambda^0_1
-\gamma M\right)<0,
$$
we get a contradiction, so that $\alpha _0=0$. It follows that
$$
u(0,0)=v(0,0)\geq \psi _0 (0,0)=\frac C{\max_{[-1,1]} \Gamma^0_1
}\left(\frac 1 2 -\frac B{2a\left(B^2+1\right)}\right)>\frac
C{3\max_{[-1,1]} \Gamma^0_1 }=:\ep _0,
$$
for  $a\geq a_0$, with $a_0>0$ sufficiently large and independent
on $b\geq 1$ and $0\leq \sigma \leq 1$. This concludes the proof
of $(iii)$.

To prove $(iv)$, observe that, since $-\mathcal R \leq \lambda
_\infty^{\delta/3}$, $\Gamma _b ^{\delta/3,\,\eta} (z):=\Gamma
_b^{\delta/3} (z)+\eta \Gamma _\infty^{\delta/3} (z)$ is a
supersolution for \eqref{eq-pour-v}. Hence we can reproduce the
proof of Lemma \ref{lem:local-pb}. Notice in particular that if
$(c_1,u_1)$ and $(c_2,u_2)$ are two solutions of $\tilde P_0(a,b)$
with $c_1>c_2$, then $u_1<u_2$.
 \fdem

\subsection{Construction of a solution in the box}\label{ss:construction-box}

Equipped with a priori estimates of Subsections \ref{ss:apriori-c}
and \ref{ss:apriori-u}, we are now in the position to construct a
solution to $P(a,b,\ep)$, with $\ep\in(0,\min(\ep^*,\ep _0))$. We
shall use a Leray-Schauder topological degree argument (see e.g.
\cite{Ber-Nic-Sch} or \cite{Ber-Nad-Per-Ryz} for related
arguments).

\begin{prop}[The solution in a box]\label{prop:sol-box} Let
$\ep\in(0,\min(\ep^*,\ep _0))$ be arbitrary. There exist $K>0$ and
$b_0>0$ such that for any $b\geq b_0$  the following holds. There
exists $a_0=a_0(b,\ep)$ such that, for all $a\geq a_0$, the
problem $P(a,b,\ep)$ has a solution $(c,u)$ such that
\begin{equation}\label{estimates-box}
\Vert u\Vert _{C^2(Q)}\leq K,\quad 0 < c\leq c^*.
\end{equation}
\end{prop}

\noindent {\bf Proof.} For a given nonnegative function $v$
defined on $Q$ and satisfying the Dirichlet boundary conditions as
requested in $P(a,b,\ep)$, consider the family $0\leq\tau\leq1$ of
linear problems
\begin{equation}\label{droite-gelee}
(P^\tau _c)\;\begin{cases}\, -\mathcal E(U)(x,z)-c U_x(x,z)\\
= \left(r(z)- \tau\displaystyle\int
_{-b}^b k(z,z')v(x,z')\,dz'-\gamma(1-\tau)v(x,z)\right)v(x,z) \quad &\text{ in }Q\vspace{8pt} \\
\,U(x,z)=\mathbf{1}_{\{x=-a\}}(x)\Gamma_b^{\delta/3}(z)\quad
&\text{ on }\partial Q.
\end{cases}
\end{equation}
Let us define $\mathcal K _\tau$ the solution operator of the
above system. More precisely $\mathcal K _\tau$ is the mapping of
the Banach space $X:=\R \times C^{1,\alpha}(Q)$ --- equipped with
the norm $\Vert (c,v)\Vert_X:=\max\left(|c|,\Vert v\Vert
_{C^{1,\alpha}}\right)$--- onto itself defined by
$$
\mathcal K _\tau:(c,v)\mapsto \left(\ep-v(0,0)+c,U^\tau _c:=\text{
the solution of } (P^\tau _c)\right).
$$
Constructing a solution of $P(a,b,\ep)$ is equivalent to showing
that the kernel of $\text{Id} -\mathcal K _1$ is nontrivial. The
operator $\mathcal K _\tau$ is compact and depends continuously on
the parameter $0\leq \tau \leq 1$. Thus the Leray-Schauder
topological argument can be applied. Define the open set
$$
S:=\left\{(c,v):\, 0< c< c^*+1,\;v>0,\; \Vert v \Vert
_{C^{1,\alpha}}< M+1\right\}\subset X,
$$
where $M>0$ is as in Lemma \ref{lem:apriori-bound}. It follows
from the a priori estimates Lemma \ref{lem:apriori-bound}, Lemma
\ref{lem:a-priori-speed} and Lemma \ref{lem:a-priori-speed2}, that
there exists $a_0=a_0(b,\ep)>0$ such that, for any $a\geq a_0$,
$0\leq \tau \leq 1$ the operator $\text{Id} -\mathcal K _\tau$
cannot vanish on the boundary $\partial S$. By the homotopy
invariance of the degree we thus have
$\text{deg}(\text{Id}-\mathcal
K_1,S,0)=\text{deg}(\text{Id}-\mathcal K_0,S,0)$. Additionally,
thanks to Lemma \ref{lem:local-pb}, any element of the kernel of
$\text{Id}-\mathcal K_0$ belongs to $S$, so that
$$
\text{deg}(\text{Id}-\mathcal
K_0,S,0)=\text{deg}(\text{Id}-\mathcal K_0,\tilde S,0),
$$
where
$$
\tilde S:=\left\{(c,v):\, -\bar c< c< c^*+1,\;v>0,\; \Vert v \Vert
_{C^{1,\alpha}}< M+1\right\}\subset X,
$$
with $\bar c =\bar c(a)>0$ as in Lemma \ref{lem:locsigma} $(iii)$.

Let us now consider the family $0\leq\sigma \leq 1$ of local and
linear problems associated with \eqref{pblocsigma}, namely
\begin{equation*}
(\tilde P^\sigma _c)\;\begin{cases}\, -\mathcal E(U)(x,z)-c
U_x(x,z)
= \left(r(z)-(1-\sigma)\mathcal R-\sigma \gamma v(x,z)\right)v(x,z) \quad &\text{ in }Q\vspace{8pt} \\
\,U(x,z)=\mathbf{1}_{\{x=-a\}}(x)\Gamma_b^{\delta/3}(z)\quad
&\text{ on }\partial Q,
\end{cases}
\end{equation*}
and let $\tilde{\mathcal K} _\sigma$ be the associated solution
operator, that is
$$
\tilde{\mathcal K} _\sigma:(c,v)\mapsto \left(\ep-v(0,0)+c,\tilde
U^\sigma _c:=\text{ the solution of } (\tilde P^\sigma _c)\right).
$$
The operator $\tilde{\mathcal K} _\sigma$ is compact and depends
continuously on the parameter $0\leq \sigma \leq 1$. The analysis
of the local problems $\tilde P_\sigma (a,b)$ in Lemma
\ref{lem:locsigma} shows that $\text{Id}-\tilde{\mathcal
K}_\sigma$ cannot vanish on the boundary of $\tilde S$. Since
${\mathcal K} _0=\tilde {\mathcal K} _1$ we have
$$
\text{deg}\left(\text{Id}-\mathcal K_0,\tilde
S,0\right)=\text{deg}\left(\text{Id}-\tilde {\mathcal K}_0,\tilde
S,0\right).
$$

To complete the proof, let us compute
$\text{deg}\left(\text{Id}-\tilde{\mathcal K}_0,\tilde S,0\right)$
by using two additional homotopies. First, consider, for $0\leq
\tau \leq 1$,
$$
\mathcal G _\tau:(c,v)\mapsto \left(\ep-(1-\tau)v(0,0)-\tau \tilde
U^0_c(0,0)+c,\tilde U^0_c:=\text{ the solution of } (\tilde P^0
_c)\right).
$$
If $\mathcal G _\tau (c,v)=(c,v)$ for some $(c,v)\in\partial
\tilde S$, then $(c,v)\in \partial \tilde S$ solves the local
problem $\tilde P^0(a,b)$ and is such that $v(0,0)=\ep$. By the a
priori estimates of Lemma \ref{lem:locsigma}, this cannot be.
Therefore $\text{Id} -\mathcal G _\tau$ does not vanish on the
boundary $\partial \tilde S$. Since $\tilde {\mathcal K}
_0=\mathcal G _0$ we have $\text{deg}\left(\text{Id}-\tilde
{\mathcal K}_0,\tilde
S,0\right)=\text{deg}\left(\text{Id}-\mathcal G_1,\tilde
S,0\right)$. Next, we know from Lemma \ref{lem:locsigma} $(iv)$
that there is a unique $(c_0,\tilde U_{c_0})\in \tilde S$ which
solves the local problem $\tilde P^0(a,b)$ and is such that
$\tilde U_{c_0 }(0,0)=\ep$. Then, consider, for $0\leq \tau \leq
1$,
$$
\mathcal H _\tau:(c,v)\mapsto \left(\ep-\tilde U^0_c(0,0)+c,\tau
\tilde U^0_c+(1-\tau) \tilde U _{c_0}\right).
$$
If $\mathcal H _\tau (c,v)=(c,v)$ for some $(c,v)\in\partial
\tilde S$, then the uniqueness in Lemma \ref{lem:locsigma} $(iv)$
 enforces $c=c_0$, $\tilde U^0 _c=\tilde U_{c_0}=v$ so that
$(c,v)\in
\partial \tilde S$ solves the local problem $\tilde P^0(a,b)$ and is such that
$v(0,0)=\ep$, which cannot be. Therefore $\text{Id} -\mathcal H
_\tau$ does not vanish on the boundary $\partial \tilde S$. Since
$\mathcal H _1=\mathcal G _1$ we have
$\text{deg}\left(\text{Id}-\mathcal G _1,\tilde
S,0\right)=\text{deg}\left(\text{Id}-\mathcal H_0 ,\tilde
S,0\right)$, where
$$
\text{Id}-\mathcal H _0:(c,v)\mapsto \left(\tilde
U^0_c(0,0)-\ep,v-\tilde U _{c_0}\right).
$$
As seen in the proof of Lemma \ref{lem:locsigma} $(iv)$, $\tilde
U^0 _c(0,0)$ is strictly decreasing in $c$ so the degree of the
first component of the above operator is $-1$. Clearly the degree
of the second one is 1. Hence $\text{deg}(\text{Id}-\mathcal
H_0,\tilde S,0)=-1$
 so that $\text{deg}(\text{Id}-\mathcal
K_1,S,0)=-1$ and there is a
 $(c,u)\in S$ solution of $P(a,b,\ep)$. This concludes the proof of
the proposition. \fdem

\section{The front with minimal speed $c^*$}\label{s:minimal-speed}

Equipped with the solution $(c,u)$ of $P(a,b,\ep)$, with
$\ep\in\left(0,\min(\ep^*,\ep _0)\right)$, of Proposition
\ref{prop:sol-box}, we now let $a\to +\infty$. Note that we have
the bounds \eqref{estimates-box} on $c$ and $u$, and also the
Gaussian control of the tails in Lemma \ref{lem:apriori-tails}.
This enables to construct
--- passing to a subsequence $a_n \to +\infty$--- a speed
$0\leq c_b \leq c^*$ and a function $u_b\in C^2 _b \left(\R\times
[-b,b]\right)$ with the same bounds as those of $u$. Similarly, we
can then consider $b\to+\infty$, to construct, via a subsequence
$b_n\to+\infty$, a speed $0\leq c \leq c^*$ and a function $u\in
C^2 _b (\R^2)$, such that $0<u\leq K$, and
\begin{equation}\label{eq-onde-construite}
-\mathcal E(u)(x,z)-c u_x(x,z)= \left(r(z)- \displaystyle\int _\R
k(z,z')u(x,z')\,dz'\right)u(x,z) \quad \text{ in } \R ^2
\end{equation}
\begin{equation}\label{masse-onde-construite}
u(0,0)=\ep
\end{equation}
\begin{equation}\label{tails-onde-construite}
0\leq u(x,z)\leq  \bar M \, \Gamma^{2\delta/3}_\infty(z),\quad
\forall (x,z)\in \R ^2.
\end{equation}

\subsection{The constructed wave has the minimal speed $c^*$}\label{ss:minimal-speed}

Here, we show that, by reducing the normalization
\eqref{masse-onde-construite} if necessary, the above constructed
solution has speed $c=c^\ast$.

\begin{lem}[A priori estimate for the infimum]\label{lem:infimum} There exists $\ep>0$ such that any solution $(c,u)$
of \eqref{eq-onde-construite}, \eqref{tails-onde-construite} with
$c\geq 0$ and $\inf_{x\in \R} u(x,0)>0$ actually satisfies $\inf
_{x\in \R} u(x,0)>\ep$.
\end{lem}

\noindent{\bf Proof.} We choose $R>0$ large enough, such that
\begin{equation}\label{truc}
\lambda_\infty^0 \leq \lambda_R^0<\frac{\lambda_\infty^0}2<0,
\quad k^+\int_{[-R,R]^c}\bar M
\Gamma^{2\delta/3}_\infty(z)\,dz\leq \frac{-\lambda_\infty^0}4.
\end{equation}
Thanks to the Harnack inequality, there exists $C>0$ such that
\begin{equation}\label{harnack}
 u(x,z')\leq C u(x,z)\quad\text{ for all } x\in\R, |z|\leq R,
 |z'|\leq R,
\end{equation}
which, combined with \eqref{tails-onde-construite} and the second
part of \eqref{truc}, implies $\int k(z,z')u(x,z')\,dz'\leq
\frac{-\lambda_\infty^0}4+2k^+CR\,u(x,z)$ in the strip $\R \times
(-R,R)$. Hence, $-\mathcal E (u)-cu_x \geq\left(
r(z)-\left(\frac{-\lambda_\infty^0}4+2k^+CR\,u\right)\right)u$ in
  $\mathbb R\times(-R,R)$. Changing variables, the function
$ v(x,y):=u\left(\frac{x-By}{\sqrt{B^2+1}},\sqrt{B^2+1}\,y\right)$
then satisfies
\begin{equation}\label{bidule}
 -v_{xx}-v_{yy}-c\sqrt{B^2+1}v_x\geq
 \left(r\left(\sqrt{B^2+1}\,y\right)-\left(\frac{-\lambda_\infty^0}4+2k^+CR\,v\right)\right)v,
\end{equation}
in $S:=\R \times \left(-\frac R {\sqrt{B^2+1}},\frac R
 {\sqrt{B^2+1}}\right)$.

Now let $\eta >0$ be arbitrarily given. Define, for $\alpha >0$,
$$
\psi _\alpha( x,y):=\alpha (1-\eta
x^2)\Gamma_R^0\left(\sqrt{B^2+1}y\right).
$$
Observe that the Harnack inequality \eqref{harnack} implies
$\psi_\alpha\leq v$ for $\alpha=C^{-1}\left\Vert \Gamma
_R^0\right\Vert _\infty ^{-1}\inf_{x\in\mathbb R}u(x,0)$, and that
\eqref{tails-onde-construite} implies  $\psi_{\alpha } (0,0)=\bar
M \geq v(0,0)$ for $\alpha=\bar M$. We can therefore define
$$
\alpha_0:=\max \left\{\alpha> 0:\,\forall (x,y)\in \bar
S,\,\psi_{\alpha}(x,y)\leq v(x,y)\right\} \in (0,\bar M].
$$
Hence $v-\psi_{\alpha _0}$ attains a zero minimum at a point
$(x_0,y_0)$ --- depending on $\eta$--- which must lie in
$\left(-\frac 1{\sqrt{\eta}},\frac 1{\sqrt{\eta}}\right) \times
\left(-\frac R{\sqrt{B^2+1}},\frac R{\sqrt{B^2+1}}\right)$ so that
$$0\geq -\Delta_{x,y}
\left(v-\psi_{\alpha_0}\right)(x_0,y_0)-c\sqrt{B^2+1} \partial _x
\left(v-\psi_{\alpha_0}\right)(x_0,y_0).$$ Using \eqref{bidule}
and straightforward computations, we arrive at
$$
0\geq \left(\frac {\lambda_\infty^0} 4-2k^+ CRv(x_0,y_0)-\lambda
_R^0\right)v(x_0,y_0)-2\alpha _0\left( \eta +c\sqrt{B^2+1} \eta
x_0\right)\Gamma _R^0\left(\sqrt{B^2+1}y_0\right).
$$
Using the first part of \eqref{truc}, $\alpha _0 \leq \bar M$ and
$|x_0|\leq \frac 1 {\sqrt \eta}$, this yields
$$
0\geq
\left(\frac{-\lambda_\infty^0}4-2k^+CRv(x_0,y_0)\right)v(x_0,y_0)-2\bar
M\left\Vert \Gamma _R^0\right\Vert _\infty\left
(\eta+c\sqrt{B^2+1}\sqrt \eta\right).
$$
It follows from the Harnack inequality \eqref{harnack} that
$v(x_0,y_0)\geq \frac 1C\inf_{x\in\mathbb R}u(x,0)>0$, so that
$$
v(x_0,y_0)\geq \frac
1{2k^+CR}\left(\frac{-\lambda_\infty^0}4-\frac{2C \bar
M\|\Gamma_R^0\|_\infty\left(\eta+c\sqrt{B^2+1}\sqrt
\eta\right)}{\inf_{x\in\mathbb R}u(x,0)}\right).
$$
Since $\eta>0$ can be chosen arbitrarily small, we have
$v(x_0,y_0)\geq \frac {-\lambda_\infty^0}{8k^+CR}$ and then
$\alpha_0\geq \ep:=\frac {-\lambda_\infty^0}{8
k^+CR\left\|\Gamma_R^0\right\|_\infty}>0$. Hence, $v(x,y)\geq \ep
(1-\eta x^2)\Gamma_R^0\left(\sqrt{B^2+1}y\right)$ for all
$(x,y)\in \bar S$. Since $\eta>0$ can be chosen arbitrarily small,
we have $v(x,y)\geq \frac{\ep} 2 \Gamma _R^0
\left(\sqrt{B^2+1}y\right)$,
 and in particular $\inf_{x\in\mathbb R}u(x,0)=\inf _{x\in \R} v(x,0)\geq \frac
 {\ep}2$. This proves the lemma. \fdem

\medskip

As a result, the constructed solution of
\eqref{eq-onde-construite}, \eqref{masse-onde-construite},
\eqref{tails-onde-construite} satisfies $\inf _{x\in \R}u(x,0)=0$.
Without loss of generality, we may assume
$\liminf_{x\to\infty}u(x,0)=0$. The following proposition then
enforces $c=c^*$ for the constructed wave. It is also of
independent interest since it proves the non existence of waves
for $0\leq c< c^*$ as stated in Theorem \ref{th:tw} $(ii)$.

\begin{prop}[$c=c^*$ for the constructed wave]\label{lem:inf-decolle} Any solution $(c,u)$
of \eqref{eq-onde-construite}, \eqref{tails-onde-construite} with
$c\geq 0$ and $\liminf_{x\to\infty}u(x,0)=0$ actually satisfies
$c\geq c^\ast$.
\end{prop}

\noindent {\bf Proof.} Assume by contradiction that $0\leq
c<c^\ast$. Choose $c<\tilde c<c^*$. Since
$\lambda_R^0\to\lambda_\infty^0=-(B^2+1)\frac{{c^\ast}^2}4$ as
$R\to\infty$, we can choose $R>0$ such that
\begin{equation}\label{inegalites}
\frac{\lambda_R^0 }{B^2+1}\leq-\frac{{\tilde c}^2}4-\frac12
\frac{{c^\ast}^2-{\tilde c}^2}4, \quad k^+\int_{[-R,R]^c} \bar M
 \Gamma^{2\delta/3}_\infty\leq\frac{B^2+1}4\frac{{c^\ast}^2-{\tilde c}^2}4.
\end{equation}Let us define the open rectangle
$$
\Omega:=\left\{(x,y):\; |x|\leq \bar x:=
\frac{\pi}{\sqrt{\left(B^2+1\right)({{\tilde c}}^2-c^2)}}\,,\,
|y|\leq \bar y:=\frac R{\sqrt{B^2+1}}\right\}.
$$
Thanks to the Harnack inequality, there exists $C>0$ such that for
any solution $(c,u)$ of \eqref{eq-onde-construite},
\eqref{tails-onde-construite} with $0\leq c \leq c^*$, and for all
$x_1 \in \R$
\begin{equation}\label{Harnack-translation}
\max _{(x,y,z')\in\bar \Omega
\times[-R,R]}u\left(x_1+\frac{x-By}{\sqrt{B^2+1}},z'\right)\leq
Cu(x_1,0).
\end{equation}

Following the change of variables of Lemma
\ref{lem:a-priori-speed}, we see that $
v(x,y):=u\left(\frac{x-By}{\sqrt{B^2+1}},\sqrt{B^2+1}\,y\right)$
satisfies
\begin{eqnarray*}
 -v_{xx}-v_{yy}-c\sqrt{B^2+1}v_x&=&\left[r\left(\sqrt{B^2+1}\,y\right)-\int
 k\left(\sqrt{B^2+1}\,y,z'\right)u\left(\frac{x-By}{\sqrt{B^2+1}},z'\right)\,dz'\right]
 v\\
 &\geq & \left[r\left(\sqrt{B^2+1}\,y\right)-\left(\frac{B^2+1}4 \frac{{c^*}^2-{\tilde c}^2}4+k^+2CRu(0,0)\right)\right]
 v,
\end{eqnarray*}
where we have used the second inequality in \eqref{inegalites} and
\eqref{Harnack-translation} with $x_1=0$. Next, define the
function
$$
\psi(x,y):=
\Gamma_R^0\left(\sqrt{B^2+1}\,y\right)e^{-\frac{c\sqrt{B^2+1}x}2}\sin\left(\frac{\sqrt{\left(B^2+1\right)({\tilde
{c}}^2-c^2)}}2x+\frac\pi 2\right).
$$
We have $\psi =0$ on $\partial \Omega$ and, using the first
inequality in \eqref{inegalites},
$$
-\psi_{xx}-\psi_{yy}-c\sqrt{B^2+1}\psi_x-r\left(\sqrt{B^2+1}\,y\right)\psi=
\left((B^2+1)\frac{{\tilde c}^2}4+\lambda_R^0\right)\psi\leq
-\frac {B^2+1}2\frac{{c^\ast}^2-{\tilde c}^2}4\psi,
$$
in $\Omega$.

Since $0<v\leq K$, we can define $\alpha
_0:=\max\{\alpha>0:\,\alpha\psi\leq v \, \text{ in } \bar
\Omega\}>0$ and there is a point $(x_0,y_0)\in \Omega$ where
$w:=\alpha _0 \psi -v$ attains a zero maximum. In view of the
above inequalities, we have at point $(x_0,y_0)$,
\begin{multline}
0\leq -w_{xx}-w_{yy}-c\sqrt{B^2+1}
w_x-r\left(\sqrt{B^2+1}\,y\right)w\\
\leq\left[-\frac {B^2+1} 4
\frac{{c^*}^2-{\tilde c }^2}4+k^+2CRu(0,0)\right]v(x_0,y_0),
\end{multline}
which in turn implies $u(0,0)\geq (B^2+1)\frac{{c^*}^2-{\tilde
c}^2}{32 k^+ CR}$. In view of \eqref{Harnack-translation}, the
argument is invariant under translation w.r.t. $x$ variable, so
that
$$
\inf_{x\in\R}u(x,0)\geq (B^2+1)\frac{{c^*}^2-{\tilde c}^2}{32 k^+
CR}>0,
$$
which contradicts $\liminf_{x\to\infty}u(x,0)=0$. The proposition
is proved.\fdem

\subsection{Behaviors as $x\to \pm \infty$}\label{ss:behaviors}

The following proposition will show that the constructed wave
satisfies the lower bound in \eqref{gauche} and \eqref{cond-bords}
in Theorem \ref{th:tw}.

\begin{prop}[Behaviors at infinity]\label{minorgauche}
Let $(c,u)$ be a solution of \eqref{eq-onde-construite},
\eqref{masse-onde-construite}, \eqref{tails-onde-construite} with
$c\geq 0$. Then the following holds.
\begin{description}
 \item $(i)$ There exist $R>0$ and $\kappa >0$ such that $u(x,z)\geq \kappa \Gamma_R^0 (z)$ for all $(x,z)\in(-\infty,0]\times
 [-R,R]$.
\item $(ii)$ If $\ep>0$ is as in Lemma \ref{lem:infimum}, then
both $\displaystyle \int _{\R} u(x,z)\,dz\to 0$ and $\max _{z\in
\R} u(x,z)\to 0$, as $x\to\infty$.
\end{description}
\end{prop}

\noindent {\bf Proof.} Let us prove $(i)$. We start as in the
proof of Lemma \ref{lem:infimum}: choose $R>0$ large enough so
that \eqref{truc} holds, choose $C>0$ such that both
\eqref{harnack} and
\begin{equation}\label{harnack1}
\min _{-\frac 1 2 \leq x\leq 0\,,\; |y|\leq \frac R{\sqrt{B^2+1}}}
v(x,y)\geq \frac 1 C u(0,0)=\frac \ep C
\end{equation}
hold, and observe that $
v(x,y):=u\left(\frac{x-By}{\sqrt{B^2+1}},\sqrt{B^2+1}\,y\right)$
satisfies \eqref{bidule} in the strip $S:=\R \times \left(-\frac R
{\sqrt{B^2+1}},\frac R
 {\sqrt{B^2+1}}\right)$. Now for $\eta >0$, we define
$$
\psi _\eta( x,y):=\alpha\left(\frac 12+\eta
x\right)\Gamma_R^0\left(\sqrt{B^2+1}\,y\right),\quad \alpha:=\min
\left(\frac{\ep }{C\left\Vert \Gamma _R^0\right\Vert
_\infty},\frac {-\lambda_\infty^0}
{5k^+CR\left\|\Gamma_R^0\right\|_\infty}\right).
$$
The definition of $\alpha$ then enforces $\psi _1 \leq v$ on
$\overline {S_-}$, where $S_-:=(-\infty,0)\times \left(-\frac R
{\sqrt{B^2+1}},\frac R
 {\sqrt{B^2+1}}\right)$. We can therefore define
$$
\eta _0:=\min\left\{\eta\geq 0:\,\forall (x,y)\in
\overline{S_-}\,,\,\psi_{\eta}(x,y)\leq v(x,y)\right\}\in [0,1].
$$
 Let us assume by contradiction that
 $\eta_ 0 >0$. Function $v-\psi_{\eta _0}$ then attains a zero minimum at a point
$(x_0,y_0)$; the definition of $\alpha$ and the Harnack inequality
\eqref{harnack1} prevents $x_0=0$ so that $(x_0,y_0)$ has to lie
in $\left(-\frac 1{2\eta_0},0\right) \times \left(-\frac
R{\sqrt{B^2+1}},\frac R{\sqrt{B^2+1}}\right)$. We therefore have
$0\geq -\Delta_{x,y} (v-\psi_{\eta_0})(x_0,y_0)-c\sqrt{B^2+1}
\partial _x (v-\psi_{\eta_0})(x_0,y_0)$. Using \eqref{bidule}, $\partial_x \psi_{\eta_0}\geq 0$ and the first part of \eqref{truc}
 we arrive at
$$
0\geq \left(\frac {-\lambda_\infty^0} 4-2k^+
CRv(x_0,y_0)\right)v(x_0,y_0),
$$
which in turn implies  $\frac{-\lambda_\infty^0}{8 k^+CR}\leq
v(x_0,y_0)=\psi_{\eta _0}(x_0,y_0)\leq \frac \alpha 2 \left\Vert
\Gamma _R^0\right\Vert _\infty$, which contradicts the definition
of $\alpha$. Hence $\eta _0=0$ and $v(x,y)\geq \frac \alpha 2
\Gamma _R^0\left(\sqrt{B^2+1}\,y\right)$ for all
$(x,y)\in(-\infty,0]\times\left[-\frac R {\sqrt{B^2+1}},\frac R
 {\sqrt{B^2+1}}\right]$.  This concludes the proof of $(i)$.

Thanks to the Harnack inequality and the control of the tails
$u(x,z)\leq \bar M \Gamma^{2\delta/3}_\infty(z)$, in order to
prove $(ii)$ it is enough to prove $u(x,0)\to$, as $x\to \infty$.
Assume by contradiction that there exists $\nu >0$ and
$x_n\to+\infty$ such that $u(x_n,0)\geq \nu$, for all $n$. Then,
the proof of $(i)$ shows that $u(x,0)\geq \frac 12\min
\left(\frac{\nu}{C\|\Gamma_R^0\|_\infty}, \frac
{-\lambda_\infty^0}{5k^+ CR\|\Gamma_R^0\|_\infty}\right)$ for all
$x\in (-\infty,x_n)$, and then for all $x\in\R$. Hence $\inf
_{x\in \R} u(x,0)>0$ so that Lemma \ref{lem:infimum} implies $\inf
_{x\in \R}u(x,0)>\ep$, which contradicts the normalization
\eqref{masse-onde-construite}. This proves $(ii)$. \fdem

\section{Faster fronts ($c>c^*$)}\label{s:faster}

In this section we fix $c>c^*= 2
\sqrt{\frac{-\lambda_\infty^0}{B^2+1}}$ and construct a
nonnegative function $u\in C^2 (\R^2)$ solution of
\begin{equation}\label{pbc}
 -\mathcal E(u)(x,z)-c u_x(x,z)= \left(r(z)- \displaystyle\int _\R
k(z,z')u(x,z')\,dz'\right)u(x,z) \quad \text{ in } \R ^2.
\end{equation}
Using the change of variables
$v(x,y):=u\left(\frac{x-By}{\sqrt{B^2+1}},\sqrt{B^2+1}\,y\right)$,
we need to construct a nonnegative $v=v(x,y)$ solution of
\begin{multline}
 \mathcal L v (x,y):=-v_{xx}(x,y)-v_{yy}(x,y)-c\sqrt{B^2+1}v_x(x,y)-r\left(\sqrt{B^2+1}\,y\right)v(x,y)=\\ -v(x,y)\int
 _{\R}k\left(\sqrt{B^2+1}\,y,z'\right)v\left(x-By+\frac
 B{\sqrt{B^2+1}}z',\frac{z'}{\sqrt{B^2+1}}\right)\,dz' \quad\text{
 in } \R^2.
\end{multline}
Note also that solving the problem in the box
\begin{equation}\label{boite}
-\mathcal E(u)(x,z)-c u_x(x,z)= \left(r(z)- \displaystyle\int
_{-b}^b k(z,z')u(x,z')\,dz'\right)u(x,z) \quad \text{ in }
Q=(-a,a)\times(-b,b),
\end{equation}
is equivalent to solving
\begin{equation}\label{eq-pour-v-boite}
 \mathcal L v (x,y)=-v(x,y)\int
 _{-b }^b k\left(\sqrt{B^2+1}\,y,z'\right)v\left(x-By+\frac
 B{\sqrt{B^2+1}}z',\frac{z'}{\sqrt{B^2+1}}\right)\,dz'\quad\text{
 in } Q_1,
\end{equation}
with $Q_1= \left\{(x,y):\,|y|< \frac
b{\sqrt{B^2+1}},\,\left|\frac{x-By}{\sqrt{B^2+1}}\right|<a\right\}$.

\medskip

We first adapt the strategy of \cite{Ber-Nad-Per-Ryz}: we
construct a solution in the box by using sub and supersolutions
and the Schauder fixed point theorem. This will allow to let $a\to
\infty$ but we shall need a extra argument to let $b\to \infty$.

\medskip

\noindent{\bf Construction of sub and supersolutions.} We use
again the supersolution of Lemma \ref{lem:a-priori-speed}: since
$c>c^*$, one can select $\mu<0$ the largest root of $\mu
^2+c\sqrt{B^2+1}\mu+\frac{{c^*}^2}4 \left(B^2+1\right)=0$. Then
the function
$$
w (x,y):=e^{\mu x}\Gamma_\infty^0\left(\sqrt{B^2+1} y\right),
$$
with $\Gamma_\infty^0$ the eigenfunction appearing in Definition
\ref{def:vp}, satisfies $\mathcal L w(x,y)=0$ in $\R^2$.

Next, we aim at constructing a kind of subsolution. Precisely we
look after a function $h$ such that (note that the supersolution
$w$ appears in the integral term)
\begin{equation}\label{sous-sol-presque}
\mathcal L h(x,y) \leq -h(x,y)\int
 _{-b}^b k\left(\sqrt{B^2+1}\,y,z'\right)w\left(x-By+\frac
 B{\sqrt{B^2+1}}z',\frac{z'}{\sqrt{B^2+1}}\right)\,dz'\quad\text{
 in } \{h>0\}.
\end{equation}
Since $k\leq k^+$ and $\int
 _{\R} e^{\mu \frac
 B{\sqrt{B^2+1}}z'} \Gamma_\infty^0 (z')\,dz'<\infty$, there is $C>0$ such
 that
\begin{equation}\label{termeint}
 \int
 _{\R}k\left(\sqrt{B^2+1}\,y,z'\right)w\left(x-By+\frac
 B{\sqrt{B^2+1}}z',\frac{z'}{\sqrt{B^2+1}}\right)\,dz'\leq C
 e^{\mu(x-By)}.
 \end{equation}
Let us choose $\ep>0$ small enough so that $-\rho :=(\mu-\ep)
^2+c\sqrt{B^2+1}(\mu-\ep)+\frac{{c^*}^2}4 \left(B^2+1\right)<0$
and $\mu+\ep<0$. For a constant $A>1$ to be selected later, let us
define
$$
h(x,y):=\left(\frac 1 A e^{\mu x}-e^{(\mu
-\ep)x}\right)\Gamma_\infty^0 \left(\sqrt{B^2+1}\,y\right)=\frac 1
A w(x,y)-e^{(\mu -\ep)x}\Gamma_\infty^0
\left(\sqrt{B^2+1}\,y\right),
$$
with is nonnegative if and only if $x\geq \ep^{-1}\ln A$. Thanks
to the estimate \eqref{termeint}, we have
\begin{align}
 &\mathcal L h (x,y)+h(x,y)\int_{-b }^b k\left(\sqrt{B^2+1}\,y,z'\right)w\left(x-By+\frac B{\sqrt{B^2+1}}z',\frac{z'}{\sqrt{B^2+1}}\right)\,dz'\nonumber\\
&\leq \left(-\rho+Ce^{\mu(x-By)}\right)h(x,y)\nonumber\\
&\leq 0,\label{presque-sous-sol}
\end{align}
in $\left\{(x,y)\in Q_1:\,x \geq \frac 1\mu\ln \left(\frac\rho
C\right)+\frac B{\sqrt{B^2+1}}b\right\}$ which contains
$\{h>0\}=\left\{(x,y)\in Q_1:\,x\geq \ep^{-1}\ln A\right\}$
provided that $A$ is chosen sufficiently large. Notice that $A$
does not depend on $a$ but does depend on $b$. We will thus need
an extra argument below.

\medskip

\noindent{\bf Construction of a propagating wave in the strip
$\mathbb R\times [-b,b]$.}  Let $b>0$ be arbitrary. Consider the
problem \eqref{eq-pour-v-boite} in the box $Q_1$ supplemented with
the boundary conditions
\begin{equation}\label{boundary}
v(x,y)=h_0(x,y):=\max\left(0,h(x,y)\right)\quad \text{ for all }
(x,y)\in
\partial Q_1.
\end{equation}
Define the convex set of functions
$$
R_{a,b}:=\left\{v\in C(Q_1):\,h_0\leq v\leq w\right\},
$$
and the compact application $\Phi_{a,b}$ that maps a given
$v^\star \in C(Q_1)$ to the solution $v$ of
$$ \mathcal L v (x,y)=-v(x,y)\int
 _{-b }^b k\left(\sqrt{B^2+1}\,y,z'\right)v^\star\left(x-By+\frac
 B{\sqrt{B^2+1}}z',\frac{z'}{\sqrt{B^2+1}}\right)\,dz'\quad\text{
 in } Q_1,
$$
supplemented with \eqref{boundary}. Since $\mathcal L v\leq
0=\mathcal L w$ in $Q_1$, $v=h_0\leq w$ on $\partial Q_1$, the
maximum principle implies $v\leq w$ in $Q_1$. Also,
\eqref{presque-sous-sol} implies $\mathcal L h_0+h_0\int k w \leq
0= \mathcal L v+v\int k v^\star \leq \mathcal L v +v \int kw$ in
$Q_1\cap \{h>0\}$, and $v\geq 0$ implies $h_0\leq v$ on $\partial
\left(Q_1\cap \{h>0\}\right)$. Hence we have $h_0\leq v$ in
$Q_1\cap\{h>0\}$ and thus in $Q_1$. Hence, $\Phi_{a,b}$ maps
$R_{a,b}$ into itself. By the Schauder fixed point theorem,
$\Phi_{a,b}$ has a fixed point $v_{a,b}\in R_{a,b}$ which solves
the problem in the box \eqref{eq-pour-v-boite} and satisfies the
boundary conditions \eqref{boundary}. Hence we are equipped with
$u_{a,b}$ solution of \eqref{boite} with
\begin{equation}\label{bords}
u_{a,b}(x,z)=h_0\left(\sqrt{B^2+1}\,x+\frac B{\sqrt{B^2+1}}z,\frac
z{\sqrt{B^2+1}}\right) \quad\text{ on } \partial Q,
\end{equation}
and
\begin{equation}\label{encadrement} h_0\left(\sqrt{B^2+1}\,x+\frac
B{\sqrt{B^2+1}}z,\frac z{\sqrt{B^2+1}}\right)\leq u_{a,b}(x,z)\leq
e^{\mu \left(\sqrt{B^2+1}\,x+\frac
B{\sqrt{B^2+1}}z\right)}\Gamma_\infty^0 (z) \quad\text{ in } \bar
Q.
\end{equation}
We claim (see proof below) that there exists $\bar M>0$ such that,
for all $a>0$, $b>0$,
\begin{equation}\label{a-prouver}
 u_{a,b}(x,z)\leq \bar M\,  \Gamma^{2\delta/3}_\infty (z)\leq \bar M \left\Vert \Gamma^{2\delta/3}_\infty \right\Vert _\infty=:M, \quad\text{
 in }  \bar Q.
\end{equation}
Now, for a given $b>0$, choose $A=A_b$ as in the construction of
$h$ above. The family $(u_{a,b})_a$ is uniformly bounded, and is
then uniformly bounded in $C^{2,\alpha}(Q)$. This allows to let
$a\to\infty$, possibly along a subsequence. In this limit, we have
$u_{a,b}\to u_b$, which is a solution of \eqref{pbc} in the strip
$S_b:=\mathbb R\times (-b,b)$, and satisfies \eqref{encadrement},
\eqref{a-prouver} in $\overline{S_b}$; in particular
\eqref{encadrement} yields
$$
u_b(x_b,0)\geq h_0\left(\frac 1 \ep \ln \left(\frac{\mu -\ep}\mu
A_b\right)\right)=:\ep _b>0,\quad x_b:= \frac 1 {\sqrt{B^2+1}}
\frac 1 \ep \ln \left(\frac {\mu -\ep}\mu A_b\right).
$$
Since the problem in the strip is invariant w.r.t. translation in
the $x$ variable, we may assume $u_b(0,0)\geq \ep_b>0$. Since
$\ep_b \to 0$ as $b\to \infty$, before letting $b\to\infty$ we
need an additional argument to get a uniform w.r.t. $b$ lower
bound for $u_b(0,0)$.

Let us now prove \eqref{a-prouver}. For $0\leq x \leq a$, it
follows from \eqref{encadrement} that the mass satisfies
$m(x):=\int _{-b}^b u_{a,b}(x,z)\,dz\leq C:=\int _{\R} e^{\mu
\frac B {\sqrt{B^2+1}} z}\Gamma_\infty^0 (z)\,dz$. For $-a\leq x
\leq 0$, observe that $u_{a,b}=h_0=0$ on $\partial
((-a,0)\times(-b,b))$ so that we can follow Lemma
\ref{lem:apriori-mass} to obtain that the mass satisfies a
Fisher-KPP inequality. Since $m(-a)=0$ and $m(0)\leq C$, the
maximum principle yields $m(x)\leq \max \left(\frac {2 \max
r}{k^-},C\right)$ for $-a\leq x \leq 0$ and thus for  $-a\leq x
\leq a$. This uniform bound for the mass enables to argue exactly
as in subsection \ref{ss:apriori-u} --- recall that $c>c^*$ has
been fixed--- to get \eqref{a-prouver}.

\medskip

\noindent{\bf Uniform lower bound for $u_b(0,0)$.} We choose $R>0$
large enough so that \eqref{truc} holds, and $C>0$ such that
\eqref{harnack1} holds. Then, for $b\geq R+1$,
$v_b(x,y):=u_b\left(\frac{x-By}{\sqrt{B^2+1}},\sqrt{B^2+1}\,y\right)$
satisfies \eqref{bidule}
in the strip $S=\R \times \left(-\frac R {\sqrt{B^2+1}},\frac R
 {\sqrt{B^2+1}}\right)$. We are therefore in the position to reproduce the proof of Proposition \ref{minorgauche} $(i)$.
Hence, there is $\kappa _b >0$, depending on $u_b(0,0)$, such that
$v_b(x,y)\geq \kappa_b \Gamma_R^0(y)$ for all $(x,y)\in
(-\infty,0]\times [-R,R]$. Since $b\geq R+1$, we can apply the
Harnack inequality to show that there exists $\ep_b>0$ such that
\begin{equation}\label{minoration}
 v_b(x,y)\geq \ep_b, \quad \forall (x,y)\in (-\infty,0]\times
 [-R,R].
\end{equation}


Let us now introduce a smooth function $\phi:\mathbb R\to[0,1]$,
such that $\phi\equiv 0$ on $(-\infty, -4]\cup[-1,\infty)$, and
$\phi\equiv 1$ on $[-3,-2]$. For $\eta >0$, $\alpha \geq 0$, let
us define
$$
\psi_{\alpha,\eta}(x,y):=\alpha \phi(\eta
x)\Gamma_{R}^0\left(\sqrt{B^2+1}y\right)\quad\textrm{ for }
(x,y)\in \bar S.
$$
In view of \eqref{a-prouver} we have $v_b\leq M$ so that, for any
$\eta> 0$, we can define
$$
\alpha_0=\alpha_0(\eta,b):=\max\left\{\alpha \geq 0:\,\forall
(x,y)\in \bar S,\,\psi_{\alpha,\eta}(x,y)\leq v_b(x,y)\right\} \in
[0,\alpha _{max}],
$$
with $\alpha _{max}:= M / \left\Vert \Gamma _R^0\right\Vert
_\infty$. Hence $v_b-\psi_{\alpha _0,\eta}$ attains a zero minimum
at a point $(x_0,y_0)$ --- depending on $\eta$ and $b$--- which
must lie in $(-\frac 4\eta,-\frac 1\eta) \times \left(-\frac
R{\sqrt{B^2+1}},\frac R{\sqrt{B^2+1}}\right)$, thanks to the
definition of $\psi_{\alpha,\eta}$ and $\phi$. Then,
$$0\geq
-\Delta_{x,y}
\left(v_b-\psi_{\alpha_0,\eta}\right)(x_0,y_0)-c\sqrt{B^2+1}
\partial _x \left(v_b-\psi_{\alpha_0,\eta}\right)(x_0,y_0).
$$
Using \eqref{bidule}, the first part of \eqref{truc} and
straightforward computations, we arrive at
\begin{eqnarray*}
& 0\geq& \left(\frac
{-\lambda_\infty^0}4-2k^+CRv_b(x_0,y_0)\right)v_b(x_0,y_0)+\\
&&\alpha _0\left(\eta ^2 \phi''(\eta x_0)+\eta
c\sqrt{B^2+1}\phi'(\eta x_0)\right)\Gamma
_R^0\left(\sqrt{B^2+1}y_0\right).
\end{eqnarray*}
Then \eqref{minoration} yields
$$
0\geq \left(\frac{-\lambda_\infty^0}4-2k^+CRv_b(x_0,y_0)-\frac
{M\left(\eta ^2 \Vert \phi ''\Vert _\infty+\eta c\sqrt{B^2+1}\Vert
\phi'\Vert _\infty\right)}{\ep_b} \right)v_b(x_0,y_0).
$$
Since $\eta>0$ can be arbitrarily small, we discover that
$v_b(x_0,y_0)\geq \frac{-\lambda_\infty^0}{8k^+CR}$, which in turn
implies $\alpha_0\geq
\ep:=\frac{-\lambda_\infty^0}{8k^+CR\|\Gamma_R^0\|_\infty}$. Since
$\psi_{\alpha _0,\eta}\left(-\frac 2 \eta ,0\right)=\alpha _0$
there is a point $x_b$ where $v_b(x_b,0)\geq \ep$, which in turn
provides a point $x_b'$ where $u_b(x_b',0)\geq \ep$. Since the
problem in the strip is invariant w.r.t. translation in the $x$
variable, we can assume without loss of generality that $x_b'=0$.
Thus we have the desired lower uniform bound $u_b(0,0)\geq \ep$.

\medskip

\noindent{\bf Conclusion.} The family $(u_{b})_b$ is uniformly
bounded, and is then uniformly bounded in
$C^{2,\alpha}\left(\mathbb R^2\right)$, and we may pass to the
limit $b\to\infty$, possibly along a subsequence.  In this limit,
we have $u_b\to u$, which is a solution of \eqref{pbc} in $\R^2$,
and satisfies \eqref{encadrement}, \eqref{a-prouver} in $\R ^2$
and $u(0,0)\geq \ep>0$. Hence, we have constructed $u$ which
satisfies \eqref{eq-dans-espace}, \eqref{controle-exp} --- which
in turn implies \eqref{cond-bords}--- and the upper bound in
\eqref{gauche}. Last, the lower bound in \eqref{gauche} follows
from Proposition \ref{minorgauche} $(i)$. This concludes the proof
of Theorem \ref{th:tw} in the case $c>c^*$. \fdem

\bigskip

\noindent \textbf{Acknowledgements.}  M. A. is supported by the
French {\it Agence Nationale de la Recherche} within the project
IDEE (ANR-2010-0112-01). G. R. is partially supported by the
French {\it Agence Nationale de la Recherche} within the project
CBDif-Fr ANR-08-BLAN-0333-01. G. R. thanks Sepideh Mirrahimi for
early discussions and computations on this problem.


\begin{thebibliography}{ABCD}

\bibitem{Alf-Cov} M. Alfaro and J. Coville, {\it Rapid travelling waves in the
nonlocal Fisher equation connect two unstable states}, Appl. Math.
Lett. {\bf 25} (2012), no. 12, 2095--2099.

\bibitem{Arnold} A. Arnold, L. Desvillettes and C. Prevost,
{\it Existence of nontrivial steady states for populations
structured with respect to space and a continuous trait}, Comm.
Pure Appl. Anal. {\bf 11} (2012), no. 1, 83--96.

\bibitem{Aro-Wei2} D. G. Aronson and H. F. Weinberger, {\it Multidimensional
nonlinear diffusion arising in population genetics}, Adv. in Math.
{\bf 30}  (1978), no. 1, 33--76.

\bibitem{Benichou} O. Benichou, V. Calvez, N. Meunier and R. Voituriez, {\it Front acceleration
by dynamic selection in Fisher population waves}, Phys. Rev. E
{\bf 86}, 041908 (2012).

\bibitem{Ber-Cha} H. Berestycki and G. Chapuisat, {\it Traveling fronts guided by the environment for
reaction-diffusion equations}, submitted.

\bibitem{berestycki-hamel-cpam}
H.~Berestycki and F.~Hamel, {\em Front propagation in periodic
excitable media}, Comm. Pure Appl. Math. {\bf 55} (2002), no. 8,
949--1032.

\bibitem{berestycki-hamel-cpam2} H. ~Berestycki and F. ~Hamel, \emph{Generalized transition waves and their properties},
Comm. Pure Appl. Math. {\bf 65}, (2012), no. 5, 592--648.

\bibitem{Berestycki-Hamel-Roques-II}
H.~Berestycki, F.~Hamel and L.~Roques, {\em Analysis of the
periodically fragmented environment model. II. Biological
invasions and pulsating travelling fronts}, J. Math. Pures Appl.
(9) {\bf 84} (2005), no. 8, 1101--1146.

\bibitem{Ber-Nad-Per-Ryz} H. Berestycki, G. Nadin, B. Perthame and
L. Ryzhik, {\it The non-local Fisher-KPP equation: travelling
waves and steady states}, Nonlinearity {\bf 22} (2009), no. 12,
2813--2844.

\bibitem{Ber-Nic-Sch} H. Berestycki, B. Nicolaenko and B.
Scheurer, {\it Traveling wave solutions to combustion models and
their singular limits}, SIAM J. Math. Anal. {\bf 16} (1985), no.
6, 1207--1242.

\bibitem{Ber-Nir3} H. Berestycki and L. Nirenberg, {\it On the method of moving planes and the
sliding method}, Bol. Soc. Brasil. Mat. (N.S.) {\bf 22} (1991),
no. 1, 1--37.

\bibitem{berestycki-nirenberg}
H.~Berestycki and  L. ~Nirenberg, \emph{Travelling fronts in
cylinders}, Ann. Inst. H. Poincar\'e Anal. Non Lin\'eaire {\bf 9}
(1992), no. 5, 497--572.

\bibitem{ber-nir-var} H. Berestycki, L. Nirenberg, and S. R. S. Varadhan, {\it The
principal eigenvalue and maximum principle for second-order
elliptic operators in general domains}, Comm. Pure Appl. Math.
{\bf 47} (1994), no. 1, 47--92.

\bibitem{Ber-Ros} H. Berestycki and L. Rossi, {\it On the principal eigenvalue of
elliptic operators in $\R ^N$ and applications}, J. Eur. Math.
Soc. {\bf  8} (2006), no. 2, 195--215.

\bibitem{Bouin} E. Bouin, V. Calvez, N. Meunier, S. Mirrahimi, B. Perthame, G. Raoul and R. Voituriez,
 {\it Invasion fronts with variable motility: phenotype selection, spatial sorting and wave acceleration}, to appear in
 C. R. Math. Acad. Sci. Paris.

\bibitem{Champagnat} N. Champagnat and S. M\'el\'eard, {\it Invasion and adaptive evolution for
 individual-based spatially structured populations}, J. Math. Biol. {\bf 55} (2007), 147--188.

\bibitem{CDM2} J.~Coville, J.~D{\'a}vila and S.~Mart{\'{\i}}nez, \emph{Nonlocal
anisotropic dispersal with monostable nonlinearity}, J.
Differential Equations {\bf 244} (2008), no. 12, 3080--3118.

\bibitem{CDM3}  J.~Coville, J.~D{\'a}vila and S.~Mart{\'{\i}}nez,
  \emph{Pulsating fronts for nonlocal dispersion and KPP nonlinearity}, Annales de l'Institut Henri Poincare (C) Non Linear Analysis,
  2012, doi:10.1016/j.anihpc.2012.07.005.

\bibitem{CD1} J. ~Coville and L. ~Dupaigne, \emph{On a non-local reaction
diffusion equation arising in population dynamics},  Proc. Roy.
Soc. Edinburgh Sect. A  {\bf 137}  (2007),  no. 4, 727--755.

\bibitem{Davis} M. B. Davis, R. G. Shaw and J. R. Etterson, {\it Evolutionary responses to changing climate}, Ecology {\bf 86} (2005), no. 7, 1704--1714.

\bibitem{Duputie} A. Duputi\'e, F. Massol, I. Chuine, M. Kirkpatrick and  O. Ronce,
 {\it How do genetic correlations affect species range shifts in a changing environment ?}, Ecol. Lett. {\bf 15} (2012), 251--259.

 \bibitem{Etterson} J. R. Etterson, D. E. Delf, T. P. Craig, Y. Ando and T. Ohgushi,
{\it Parallel patterns of clinal variation in Solidago altissima
in its native range in central U.S.A. and its invasive range in
Japan}, Botany {\bf 86} (2007), 91--97.

\bibitem{Fis} R. A. Fisher, {\it The wave of advance of advantageous genes},
Ann. of Eugenics {\bf 7} (1937), 355--369.

\bibitem{Gen-Vol-Aug} S. Genieys, V. Volpert and P. Auger, {\it
Pattern and waves for a model in population dynamics with nonlocal
consumption of resources}, Math. Model. Nat. Phenom. {\bf 1}
(2006), no. 1, 65--82.

\bibitem{Gil-Tru} D. Gilbarg and N. Trudinger, {\it Elliptic Partial
Differential Equations of Second Order}, Springer-Verlag: Berlin,
1977.

\bibitem{Griffith} T. M. Griffith and M. A. Watson, {\it Is evolution necessary for range expansion? Manipulating
reproductive timing of a weedy annual transplanted beyond its
range}, Am. Nat. {\bf 167} (2006), no. 2, 153--164.

\bibitem{Hermsen} R. Hermsen, J. B. Deris and T. Hwa,
 {\it On the rapidity of antibiotic resistance evolution facilitated by a concentration gradient}, Proc. Nat. Acad. Sci. USA {\bf 109} (2012), 10775--10780.

\bibitem{HZ} W. Hudson and B. Zinner, \emph{Existence of traveling waves for
reaction diffusion equations of {F}isher type in periodic media},
Boundary value problems for functional-differential equations,
187--199, World Sci. Publ., River Edge, NJ, 1995.

\bibitem{Jabin} P. E. Jabin and G. Raoul, {\it On selection dynamics for competitive interactions}, J. Math. Biol. {\bf 63} (2011), no.
3, 493--517.

\bibitem{Keller} S. R. Keller and D. R. Taylor, {\it History, chance and adaptation during biological
invasion: separating stochastic phenotypic evolution from response
to selection}, Ecol. lett. {\bf 11} (2008), 852--866.

\bibitem{Keymer} J. E. Keymer, P. Galajda, C. Muldoon, S. Park and R. H. Austin,
 {\it Bacterial metapopulations in nanofabricated landscapes}, Proc. Nat. Acad. Sci. USA {\bf 103} (2006), no. 46, 17290--17295.

\bibitem{Kirkpatrick} M. Kirkpatrick and N. H. Barton, {\it Evolution of a species' range}, Amer. Nat. {\bf 150} (1997), no. 1, 1--23.

\bibitem{KPP} A. N. Kolmogorov, I. G.  Petrovsky and N. S.  Piskunov, \emph{\'Etude de l'\'equation de la diffusion avec croissance de
la quantit\'e de mati\`ere et son application \`a un probl\`eme
biologique}, Bulletin Universit\'e d'\'Etat \`a Moscow (Bjul.
Moskowskogo Gos. Univ), S\'erie Internationale, (1937), Section A,
1--26.

\bibitem{Lorz} A. Lorz, S. Mirrahimi and B. Perthame, {\it Dirac mass dynamics in a multidimensional nonlocal parabolic equation}, Comm. Partial Differential
Equations {\bf 36} (2011), no. 6, 1071--1098.

\bibitem{Mir-Rao} S. Mirrahimi and G. Raoul, {\it Population
structured by a space variable and a phenotypical trait},
submitted.

\bibitem{nadin} G. Nadin, {\em Traveling fronts in space-time periodic media}, J.
Math. Pures Appl. (9) {\bf 92} (2009), no. 3, 232--262.

\bibitem{Nolen-Ryzhik} J. Nolen and L. Ryzhik, {\em Traveling waves in a one-dimensional
 heterogeneous medium}, Ann. Inst. H. Poincar\'e Anal. Non Lin\'eaire {\bf 26} (2009), no. 3, 1021--1047.

\bibitem{Peck} J. R. Peck, J. M. Yearsley and D. Waxman, {\it Explaining the geographic distributions of sexual and asexual populations}, Nature {\bf 391} (1998), 889--892.

\bibitem{Phillips} B. L. Phillips, G. P. Brown, J. K. Webb and R. Shine, {\it Invasion and the evolution of speed in toads}, Nature {\bf
439} (803) (2006).

\bibitem{Polechova} J. Polechov\'a and N. Barton, {\it Speciation through competition: a critical review}, Evolution {\bf 59} (2005), 1194--1210.

\bibitem{Prevost} C. Prevost, {\it Applications of partial differential equations and their numerical simulations of population dynamics}, PhD Thesis, University of Orleans (2004).


\bibitem{Shen2012} W. Shen and A. Zhang, {\it Traveling wave
solutions of spatially periodic nonlocal monostable equations},
ArXiv e-prints, (2012). http://arxiv.org/abs/1202.2452

\bibitem{W1} H.F. Weinberger, {\em Long-time behavior of a class of biological
models}, SIAM J. Math. Anal. {\bf 13} (1982), no. 3, 353--396.

\bibitem{Xin3} J. Xin, \emph{Front propagation in heterogeneous media}, SIAM Rev. {\bf 42}, (2000), no. 2, 161--230.


\end{thebibliography}
\end{document}